\documentclass[12pt,twoside]{article}
\usepackage[english]{babel}
\usepackage{amsmath}
\usepackage{times,amssymb,amscd}
\usepackage[utf8]{inputenc}
\usepackage[pass,legalpaper]{geometry}
\usepackage{fancyhdr}
\usepackage{multicol}
\usepackage{enumitem}
\usepackage{graphicx}
\usepackage{pgf}
\usepackage{textcomp}
\usepackage{stackengine}
\newcommand\xrowht[2][0]{\addstackgap[.5\dimexpr#2\relax]{\vphantom{#1}}}
\newtheorem{thm}{Theorem}[section]

\newtheorem{lem}[thm]{Lemma}
\newtheorem{prop}[thm]{Proposition}
\newtheorem{defn}[thm]{Definition}
\newtheorem{rem}[thm]{Remark}
\numberwithin{equation}{section}

\newcommand{\bG}{\mathbf{G}}

\newcommand{\ba}{\mathbf{a}}
\newcommand{\bc}{\mathbf{c}}

\newcommand{\bh}{\mathbf{h}}

\newcommand{\bb}{\mathbf{b}}

\newcommand{\bbH}{\mathbb{H}}

\newcommand{\Ba}{\boldsymbol{a}}
\newcommand{\Bb}{\boldsymbol{b}}

\newcommand{\Bu}{\boldsymbol{u}}

\newcommand{\cD}{\mathcal{D}}
\newcommand{\cF}{\mathcal{F}}

\newcommand{\cT}{\mathcal{T}}
\newcommand{\cC}{\mathcal{C}}

\newcommand{\cB}{\mathcal{B}}

\begin{document}
\pagestyle{myheadings}
\markboth{\centerline{Arnasli Yahya and Jen\H o Szirmai}}
{Optimal ball and horoball packings...}
\title
{Optimal Ball and Horoball Packings Generated
	by Simply Truncated Coxeter
	Orthoschemes with Parallel Faces in Hyperbolic $n$-space for $4 \leq n \leq 6$
 }
\author{Arnasli Yahya and Jen\H o Szirmai \\
\normalsize Department of Geometry, Institute of Mathematics,\\
\normalsize Budapest University of Technology and Economics, \\
\normalsize M\H uegyetem rkp. 3., H-1111 Budapest, Hungary \\
\normalsize arnasli@math.bme.hu,~szirmai@math.bme.hu
\date{\normalsize{\today}}}
\maketitle
\begin{abstract} 
   After investigating the $3$-dimensional case \cite{YSz22}, we continue to address and close the problems of optimal ball and horoball packings in truncated Coxeter orthoschemes with parallel faces that exist in $n$-dimensional hyperbolic space $\overline{\mathbb{H}}^n$ up to $n=6$. In this paper, we determine the optimal ball and horoball packing configurations and their densities for the aforementioned tilings in dimensions $4 \leq n \leq 6$, where the symmetries of the packings are derived from the considered Coxeter orthoscheme groups. Moreover, for each optimal horoball packing, we determine the parameter related to the corresponding Busemann function, which provides an isometry-invariant description of different optimal horoball packing configurations. 
\end{abstract}
\section{Background and preliminaries} \label{section1}
According to the known results in $n$-dimensional hyperbolic space $\bbH^n$ ($2 \leq n \in \mathbb{N}$), the densest packings of spheres (balls), horospheres (horoballs), and hyperspheres (hyperballs) are achieved through Coxeter simplex tilings. These findings provide a strong motivation for further research on these types of tilings and their associated packings.
To begin, let's review the previous results and background of the topic.

In the realm of discrete geometry, let $X$ represent a space of constant curvature, which can be either the $n$-dimensional sphere $\mathbb{S}^n$, Euclidean space $\mathbb{E}^n$, or hyperbolic space $\mathbb{H}^n$ with $n \geq 2$. One fundamental question in this field is to determine the maximum packing density in $X$ achieved by congruent, non-overlapping balls of a given radius \cite{G--K}, \cite{Fejes2}.
Defining the packing density in hyperbolic space poses certain challenges, as illustrated by B\"or\"oczky's work \cite{B78} and the notable constructions described in references \cite{G--K} or \cite{R06}. To address these challenges, the most widely accepted notion of packing density takes into account the local densities of the balls relative to their Dirichlet-Voronoi cells (see \cite{B78} and \cite{K98}). In the study of ball and horoball packings in $\overline{\mathbb{H}}^n$, we employ an extended concept of such local density.

Let $B$ be a horoball belonging to a packing $\mathcal{B}$, and let $P \in \overline{\mathbb{H}}^n$ be an arbitrary point. We define $d(P,B)$ as the shortest distance from point $P$ to the horosphere $S = \partial B$. It should be noted that if $P$ lies within the horoball $B$, then $d(P,B)\leq 0$.
The Dirichlet-Voronoi cell $\mathcal{D}(B,\mathcal{B})$ associated with the horoball $B$ can be understood as the convex body or region consisting of all points in $\overline{\mathbb{H}}^n$ whose distance to $B$ is smaller than their distance to any other horoball in the packing.
\begin{equation}
	\cD(B,\cB) = \{ P \in \mathbb{H}^n | d(P,B) \le d(P,B'), ~ \forall B' \in \cB \}. \notag
\end{equation}
Both $B$ and $\cD$ have infinite volume, therefore the standard notion of local density should be 
modified. Let $Q \in \partial{\mathbb{H}}^n$ denote the ideal center of $B$, and take its boundary $S$ to be the one-point compactification of Euclidean $(n-1)$-space.
Let $B_C^{n-1}(r) \subset S$ be the Euclidean $(n-1)$-ball with center $C \in S \setminus \{Q\}$.
Then $Q$ and $B_C^{n-1}(r)$ determine a convex cone 
$\cC^n(r) = cone_Q\left(B_C^{n-1}(r)\right) \in \overline{\mathbb{H}}^n$ with
apex $Q$ consisting of all hyperbolic geodesics passing through $B_C^{n-1}(r)$ with limit point $Q$. The local density $\delta_n(B, \cB)$ of $B$ to $\cD$ is defined as
\begin{equation}
	\delta_n(\cB, B) =\varlimsup\limits_{r \rightarrow \infty} \frac{\mathrm{vol}(B \cap \cC^n(r))} {\mathrm{vol}(\cD \cap \cC^n(r))}. \notag
\end{equation}
This limit is independent of the choice of center $C$ for $B^{n-1}_C(r)$.

In the case of periodic ball or horoball packings, this local density defined above extends to the entire hyperbolic space via its symmetry group, and 
is related to the simplicial density function (defined below) that we generalized in \cite{Sz12} and \cite{Sz12-2}.
In this paper, we shall use such definition of packing density.

A Coxeter simplex is a top dimensional simplex in $X$ with dihedral angles either integral submultiples of $\pi$ or $0$. 
The group generated by reflections on the sides of a Coxeter simplex is a Coxeter simplex reflection group. 
Such reflections generate a discrete group of isometries of $X$ with the Coxeter simplex as the fundamental domain; 
hence the groups give regular tessellations of $X$, if the fundamental simplex is characteristic. 
The Coxeter groups are finite for $\mathbb{S}^n$, and infinite for $\mathbb{E}^n$ or $\overline{\mathbb{H}}^n$.
There are non-compact Coxeter simplices in $\overline{\mathbb{H}}^n$ with ideal vertices on $\partial \mathbb{H}^n$, however only for dimensions $2 \leq n \leq 9$; 
furthermore, only a finite number exists in dimensions $n \geq 3$. 
Johnson {\it et al.} \cite{JKRT} found the volumes of all Coxeter simplices in hyperbolic $n$-space. 
Such simplices are the most elementary building blocks of hyperbolic manifolds,
the volume of which is an important topological invariant. Of course, if we allow the simplexes to have ``outer" vertices, 
the number of possible simplexes expands and we get nice new tilings (see e.g. \cite{IH1,IH2}). Our work is also related to one of these classes.

In $n$-dimensional space $X$ of constant curvature
$(n \geq 2)$, define the simplicial density function $d_n(r)$ to be the density of $n+1$ mutually tangent balls of radius $r$ in the simplex spanned by their centers. 
L.~Fejes T\'oth and H.~S.~M.~Coxeter
conjectured that the packing density of balls of radius $r$ in $X$ cannot exceed $d_n(r)$.
Rogers \cite{Ro64} proved this conjecture in Euclidean space $\mathbb{E}^n$.
The $2$-dimensional spherical case was settled by L.~Fejes T\'oth \cite{FTL}, and B\"or\"oczky \cite{B78} gave a proof for the extension.
In hyperbolic 3-space, 
the monotonicity of $d_3(r)$ was proved by B\"or\"oczky and Florian
in \cite{B--F64}; in \cite{Ma99} Marshall 
showed that for sufficiently large $n$, 
function $d_n(r)$ is strictly increasing in variable $r$. Kellerhals \cite{K98} showed $d_n(r)<d_{n-1}(r)$, and that in cases considered by 
Marshall the local density of each ball in its Dirichlet--Voronoi cell is bounded above by the simplicial horoball density $d_n(\infty)$. 

The simplicial ball and horoball packing density upper bound 
$d_3(\infty) = (1+\frac{1}{2^2}-\frac{1}{4^2}-\frac{1}{5^2}+\frac{1}{7^2}+\frac{1}{8^2}-\frac{1}{10^2}-\frac{1}{11^2}+\dots)^{-1} = 0.85327\dots$ cannot be achieved 
by packing regular balls, instead, it is realized by horoball packings of
$\overline{\mathbb{H}}^3$, the regular ideal simplex tiles $\overline{\mathbb{H}}^3$.
More precisely, the centers of horoballs in  $\partial\overline{\mathbb{H}}^3$ lie at the vertices of the ideal regular Coxeter 
simplex tiling with Schl\"afli symbol $\{3,3,6\}$. 

In three dimensions, the B\"or\"oczky-type bound for horoball packings are used for volume estimates of cusped hyperbolic manifolds 
\cite{A87,M86}, more recently \cite{ACS,MaM}. Lifts of horoball neighborhoods of cusps give horoball packings in the universal 
cover $\bbH^n$, and for some discrete torsion-free subgroup of isometries $\mathbb{H}^n/\Gamma$ is a cusped hyperbolic manifold 
where the cusps lift to ideal vertices of the fundamental domain. In this setting, a manifold with a single cusp has a well-defined maximal cusp neighborhood, while manifolds with multiple cusps have a range of non-overlapping cusp neighborhoods with 
boundaries with nonempty tangential intersection, these lift to different horoball types in the universal cover. 
An important application is Adams' proof that the Geiseking manifold is the noncompact hyperbolic $3$-manifold 
of minimal volume \cite{A87} (further information \cite{MPSz21}). Kellerhals then used the B\"or\"oczky-type bounds to estimate volumes of higher dimensional hyperbolic manifolds \cite{K98_2}.  

In \cite{KSz}, we proved that the classical horoball packing configuration in $\overline{\mathbb{H}}^3$ realizing 
the B\"or\"oczky-type upper bound  is not unique. We gave several examples of different regular horoball 
packings using horoballs of different types, that is horoballs that have different relative densities with respect 
to the fundamental domain, that yield the B\"or\"oczky--Florian-type simplicial upper bound \cite{B--F64}.

Furthermore, in \cite{Sz12,Sz12-2} we found that 
by allowing horoballs of different types at each vertex of a totally asymptotic simplex and generalizing 
the simplicial density function to $\overline{\mathbb{H}}^n$ for $n \ge 2$,
the B\"or\"oczky-type density 
upper bound is not valid for the fully asymptotic simplices for $n \geq 4$. 
In $\overline{\mathbb{H}}^4$, the locally optimal simplicial packing density is $0.77038\dots$, 
higher than the B\"or\"oczky-type density upper bound of $d_4(\infty) = 0.73046\dots$ using horoballs of a single type. 
However, these ball packing configurations are only locally optimal and cannot be extended to the entirety 
of $\overline{\mathbb{H}}^n$. Further interesting constructions are described in \cite{Sz05-2,Sz07-1}.

In a series of $5$ papers, the optimal horoball packing configurations and their densities of Koszul-type noncompact Coxeter simplex tilings were 
investigated by R.~T.~Kozma and the second author where these tilings exist in $\overline{\mathbb{H}}^n$ for $2 \le n \le 9$ (see \cite{KSz,KSz14,KSz18,KSz23}). 
In \cite{KSz23-1}, we studied some further new aspects of the topic.

We would highlight only one of them here, in \cite{KSz14}, we found seven horoball packings of Coxeter simplex tilings in $\overline{\mathbb{H}}^4$ 
that yield densities of $5\sqrt{2}/\pi^2 \approx 0.71645$, counterexamples to L. Fejes T\'oth's conjecture for the maximal packing density of $\frac{5-\sqrt{5}}{4} \approx 0.69098$ 
in his foundational book {\it Regular Figures} \cite[p. 323]{FTL}. We note here, that this density is realized also 
in our case with Schl\"afli symbol $\{4,4,3,4,\infty\}$ (see Theorem 4.4). 

In Table 1, we summarize the Böröczky-type upper bound densities $d_n(\infty)$ (see \cite{K98}) and the densities of the densest known horoball 
arrangements which we found in our previous papers.
\begin{table}[h!]
	\footnotesize{\begin{tabular}{l|l|l|l|l}
			\hline
			$n$ & Optimal Coxeter simplex packing density & Numerical Value & $d_n(\infty)$ & $\Delta$ \\
			\hline
			3 & $\left( 1+\frac{1}{2^2} - \frac{1}{4^2} - \frac{1}{5^2} + \frac{1}{7^2} + \frac{1}{8^2} - \dots \right)^{-1}$ & 0.85328\dots & 0.85328\dots & 0\\
			4 & $5 \sqrt{2} /\pi ^{2} $  & 0.71644\dots & 0.73046\dots & 0.0140\dots\\
			5 & $\displaystyle 5 / \left(7 \zeta(3)\right)$ & 0.59421\dots & 0.60695\dots & 0.0127\dots \\
			6 & $81/ \left(4 \sqrt{2} \pi ^{3} \right)$ & 0.46180\dots & 0.49339\dots & 0.0315\dots \\
			7 & $\displaystyle 28/ \left(81 \text{L}(4,3)\right)$ & 0.36773\dots & 0.39441\dots & 0.0266\dots\\
			8 & $\displaystyle 225/\left(8\pi^4\right)$ & 0.288731\dots & 0.31114\dots & 0.0223\dots\\
			9 & $1/\left( 4 \zeta(5) \right) $ & 0.24109\dots & 0.24285\dots & 0.0017\dots\\
			\hline
	\end{tabular}}
	\caption{Packing density upper and lower bounds for $\overline{\bbH}^n$, where 
		$\zeta(\cdot)$ is the Riemann Zeta function, 
		$\text{L}(\cdot, \cdot)$ is the Dirichlet $L$-Series and $\Delta$ is the gap between the packing density upper bound and our effective lower bounds.}
	\label{table:summary}
\end{table}

In \cite{YSz22}, we considered the ball and horoball packings belonging to $3$-dimensional Coxeter tilings that are 
derived by simply truncated orthoschemes with parallel faces. 
We determined the optimal ball and horoball packing arrangements and their densities
for all above Coxeter tilings in hyperbolic $3$-space $\mathbb{H}^3$.
The centers of horoballs are required to lie at ideal vertices of the
polyhedral cells constituting the tiling, and we allow horoballs of
different types at the various vertices (using Busemann function parametrization of horoballs).  
We proved that the densest packings are realized by horoballs related to 
tiling $\{\infty,3,6,\infty \}$ and to its dual $\{\infty,6,3,\infty \}$ with density $\approx 0.8413392$.

{\it In this paper, we continue and close 
	the problem of optimal ball and horoball packings related to truncated Coxeter orthoschemes 
	with parallel faces that exist in $n$-dimensional hyperbolic space $\overline{\mathbb{H}}^n$ 
	up to $n=6$, determine the optimal ball and horoball packing configurations and their densities of the above 
	tilings in dimensions $4 \leq n \leq 6$. We prove, that the optimal ball and horoball densities coincide with the known density lower bounds (see Table 1) 
	in $4$- and $5$ dimensions and 
	less in $6$-dimension. Particularly in $\overline{\mathbb{H}}^4$, it is attained 
	by two horoballs related to $\{4,4,3,4, \infty\}$ tiling with density $\frac{5\sqrt{2}}{\pi^2}\approx 0.71645$. 
	In $\overline{\mathbb{H}}^5$, $\frac{5}{7\zeta{3}}\approx 0.59421$, realized by horoballs related to $\{4,3,4,3,3,\infty\}$ Coxeter orthoscheme tiling. 
	While in $\overline{\mathbb{H}}^6$, we have optimum value $\frac{15\sqrt{2}+18}{4\pi^3}\approx 0.31617$, 
	reached by horoball packing related to Coxeter tiling $\{3,4,3,3,3,4,\infty\}$.
	
	Our method for computing densities in the projective Beltrami--Cayley--Klein model is mostly similar to the earlier $3$-dimensional case (see \cite{YSz22, ASz}), 
	although the computations are much more complicated and require new methods, especially in volume calculations.}

\section{Basic Notions}
\subsection{The projective model of hyperbolic space $\mathbb{H}^n$}
Let $\mathbb{E}^{1,n}$ denote $\mathbb{R}^{n+1}$ with the Lorentzian inner product 
$\langle \mathbf{x}, \mathbf{y} \rangle = -x^0y^0+x^1y^1+ \dots + x^n y^n \label{bilinear_form}$
where non-zero real vectors 
$\mathbf{x, y}\in\mathbb{R}^{n+1}$ represent points in projective space 
$\mathbb{P}^n=\mathbb{P}(\mathbb{E}^{n+1})$, equipped with the quotient topology of the natural projection $\Pi: \mathbb{E}^{n+1}\setminus \{\mathbf{0}\} \rightarrow \mathbb{P}^n$. Partitioning $\mathbb{E}^{1,n}$ into $Q_+=\{\mathbf v \in \mathbb{R}^{n+1}|\langle \mathbf v, \mathbf v \rangle >0\}$, $Q_0=\{\mathbf v|\langle \mathbf v, \mathbf v \rangle =0\}$, and $Q_-=\{\mathbf v|\langle \mathbf v, \mathbf v \rangle <0\}$, the proper points of hyperbolic $n$-space are $\mathbb{H}^n = \Pi(Q_-)$, $\partial \mathbb{H}^n = \Pi(Q_0)$ are the boundary or ideal points, we will refer to points in $\Pi(Q_+)$ as outer points, and $ \overline{\mathbb{H}}^n = \mathbb{H}^n \cup \partial \mathbb{H}^n$ as extended hyperbolic space.  

Points $[\mathbf{x}], [\mathbf{y}] \in \mathbb{P}^n$ are conjugate when $\langle
\mathbf{x}, \mathbf{y} \rangle = 0$. The set of all points conjugate
to $[\mathbf{x}]$ form a projective (polar) hyperplane
$Pol([\mathbf{x}])=\{[\mathbf{y}] \in\mathbb{P}^n | \langle  \mathbf{x}, \mathbf{y} \rangle =0 \}.$
Hence $Q_0$ induces a duality $\mathbb{R}^{n+1} \leftrightarrow
\mathbb{R}_{n+1}$
between the points and hyperplanes of $\mathbb{P}^n$.
Point $[\mathbf{x}]$ and hyperplane $[\Ba]$ are incident if the value of
the linear form $\Ba$ evaluated on vector $\mathbf{x}$ is
zero, i.e. $\mathbf{x}\Ba=0= \langle \mathbf{x}, \mathbf{a} \rangle$ where $\mathbf{x} \in \
\mathbb{R}^{n+1} \setminus \{\mathbf{0}\}$, and $\Ba \in
\mathbb{R}_{n
	+1} \setminus \{\mathbf{0}\}$.
Similarly, the lines in $\mathbb{P}^n$ are given by
2-subspaces of $\mathbb{R}^{n+1}$ or dual $(n-1)$-subspaces of $\mathbb{R}_{n+1}$ \cite{Mol97}.

In this paper, we set the sectional curvature of $\mathbb{H}^n$,
$K=-k^2$, to be $k=1$. The distance $d$ between two proper points
$[\mathbf{x}]$ and $[\mathbf{y}]$ is given by
\begin{equation}
	\cosh{{d}}=\frac{-\langle  \mathbf{x},\mathbf{y} \rangle }{\sqrt{\langle  \mathbf{x},\mathbf{x} \rangle
			\langle \mathbf{y},\mathbf{y} \rangle }}.
	\label{prop_dist}
\end{equation}
The projection $[\mathbf{y}]$ of point $[\mathbf{x}]$ on plane $[\mathbf{u}]$ is given by
\begin{equation}
	\mathbf{y} = \mathbf{x} - \frac{ \langle \mathbf{x}, \mathbf{u} \rangle }{\langle \Bu, \Bu \rangle} \mathbf{u},
	\label{perp_foot}
\end{equation}
where $[\mathbf{u}]$ is the pole of the hypeplane $[\Bu]$.
\subsection{On Horoballs}
A horosphere in $\overline{\mathbb{H}}^n$ ($n \ge 2)$ is as 
hyperbolic $n$-sphere with infinite radius centered 
at an ideal point $\xi \in \partial \mathbb{H}^n$ obtained as a limit of spheres through $x \in \bbH^n$ as its center $c \rightarrow \xi$. Equivalently, a horosphere is an $(n-1)$-surface orthogonal to
the set of parallel straight lines passing through $\xi \in \partial \mathbb{H}^n$. 
A horoball is a horosphere together with its interior. 

To derive the equation of a horosphere, we fix a projective 
coordinate system for $\mathbb{P}^n$ with standard basis 
$\bold{a}_i, 0 \leq i \leq n $ so that the Cayley--Klein ball model of $\overline{\mathbb{H}}^n$ 
is centered at $O = (1,0,0,\dots, 0)$, and orient it by setting point $\xi \in \partial \bbH^n$ to lie at $A_0=(1,0,\dots, 0,1)$. 
The equation of a horosphere with center
$\xi = A_0$ passing through interior point $S=(1,0,\dots,0,s)$ is derived from the equation of the 
the boundary sphere $-x^0 x^0 +x^1 x^1+x^2 x^2+\dots + x^n x^n = 0$, and the plane $x^0-x^n=0$ tangent to the boundary sphere at $\xi = A_0$. 
The general equation of the horosphere is
\begin{equation}
	0=\lambda (-x^0 x^0 +x^1 x^1+x^2 x^2+\dots + x^n x^n)+\mu{(x^0-x^n)}^2.
	\label{horopshere_eq}
\end{equation}
Evaluating at $S$ obtain
\begin{equation}
	\lambda (-1+s^2)+\mu {(-1+s)}^2=0 \text{~~and~~} \frac{\lambda}{\mu}=\frac{1-s}{1+s}. \notag
\end{equation}
For $s \neq \pm1$, the equation of a horosphere in projective coordinates is
\begin{align}
	\label{eq:horosphere}
	(s-1)\left(-x^0 x^0 +\sum_{i=1}^n (x^i)^2\right)-(1+s){(x^0-x^n)}^2 & =0.
\end{align}
\begin{rem}
The last form can be considered in the spherical coordinate parameterization 
(It is powerful to provide a horosphere visualization.) 
\begin{equation}
\begin{aligned}
&h_1=\sqrt{\frac{1-s}{2}}\sin{\theta_1}\sin{\theta_2} \cdots \sin{\theta_{n-1}}\\
&h_2=\sqrt{\frac{1-s}{2}}\sin{\theta_1}\sin{\theta_2} \cdots \cos{\theta_{n-1}}\\
&\vdots  \hspace{2 cm}\vdots\\
&h_{n-1}=\sqrt{\frac{1-s}{2}}\sin{\theta_1}\cos{\theta_2}\\
&h_n=\frac{1+s}{2}+\frac{1-s}{2}\cos{\theta_1},
\end{aligned} \notag
\end{equation}
where $ 0\le \theta_i \le \pi, i \in \{1 \dots n-2\} $ and $0 \le \theta_{n-1} < 2\pi$. 
\end{rem}
In $\overline{\mathbb{H}}^n$, there exists an isometry $g \in \text{Isom}(\bbH^n)$ for any two horoballs $B$ and $B'$  such that $g.B = B'$.
However, it is often useful to distinguish between certain horoballs of a packing; we shall use the notion of horoball type with respect to the fundamental domain of a tiling (lattice) as introduced in \cite{Sz12-2}. In \cite{KSz23}, we showed that this coincides with the Busemann function up to scaling, hence is isometry invariant. 

Two horoballs of a horoball packing are said to be of the {\it same type} or {\it equipacked} if 
and only if their local packing densities with respect to a particular cell (in our case a 
Coxeter simplex) are equal, otherwise, the two horoballs are of {\it different type}. 
For example, the horoballs centered at $A_0$ passing through $S$ with different values for the final coordinate $s \in (-1,1)$ are of different type relative to 
a given cell and the set of all horoballs centered at vertex $A_0$ is a one-parameter family.

Volumes of horoball pieces are given by J\'anos Bolyai's classical formulas.
The hyperbolic length $L(x)$ of a horospherical arc corresponded to a chord segment of length $x$ is
\begin{equation}
	\label{eq:horo_dist}
	L(x)=2 \sinh{\left(\tfrac{x}{2}\right)} .
\end{equation}
As the intrinsic geometry of a horosphere is Euclidean, the $(n-1)$-dimensional volume $\mathrm{vol}\mathcal{A}$ of a polyhedron $\mathcal{A}$ 
on the 
surface of the horosphere can be calculated as in $\mathbb{E}^{n-1}$.
The volume of the horoball piece $\mathcal{H}(\mathcal{A})$ bounded by $\mathcal{A}$, and the set of geodesic segments
joining $\mathcal{A}$ to the center of the horoball, is
\begin{equation}
	\label{eq:bolyai}
	\mathrm{vol}(\mathcal{H}(\mathcal{A})) = \tfrac{1}{n-1}\mathrm{vol}\mathcal{A}.
\end{equation}
\subsection{Coxeter orthoschemes and tilings}
In the $\overline{\mathbb{H}}^n$ $(2\leq n\in \mathbb{N})$ space, a complete orthoscheme 
$\mathcal{S}$ of degree $d$ $(0\leq d\leq2)$ is a polytope bounded 
with hyperplanes $H^0,...,H^{n+d}$, for which $H^i\perp H^j$, if $j\neq i-1,i,i+1$.
We denote $A_i$ as the vertex opposite to the hyperplane $H^i$ where ($0\leq i\leq n$).
Equivalently, the orthoscheme can be described by the sequence of the vertices of the orthoschemes 
$A_0,...,A_n$, where $A_{i-1}A_{i}$ edge is perpendicular 
to $A_{i}A_{i+1}$ edge for all $i\in\{1,...,n-1\}$. 
Here, $A_0$ and $A_n$ are called the principal vertices of the 
orthoschemes.  

Now, consider the reflections on the facets of the simply truncated orthoscheme $H^1, \cdot, H^{n+1}$, 
and denote them with  $r_0,...,r_{n+1}$, hence define the group 
\begin{equation*}
\mathbf{G}=\langle r_0,...,r_{n+1} | (r_ir_j)^{m_{ij}}=1 \rangle,
\end{equation*}
where $\alpha^{ij}=\frac{\pi}{m_{ij}}$, and $m_{ii}=1$. Note that, if $m_{ij}=\infty$
(i.e. $H^i$ and $H^j$ are parallel), then the $r_i,r_j$ pair belongs to no relation. Here, $2\leq m_{ij}\in\{\mathbb{N}\cup\infty\}$, $i\neq j$. 
Since the Coxeter group $\bG$ acts on hyperbolic 
space $\overline{\mathbb{H}}^n$ properly discontinuously, 
then the images of the orthoscheme under this action provide a 
$\mathcal{T}$ tiling of $\overline{\mathbb{H}}^n$ 
i.e. the images of the orthoscheme fills $\overline{\mathbb{H}}^n$ without overlap.

Let $S \subset \overline{\mathbb{H}}^n$ denote a complete orthoscheme bounded by
a finite set of hyperplanes $H^i$ with unit normal forms
$\mbox{\boldmath$b$}^i \in \mbox{\boldmath$V$}\!_{n+1}$ directed
outwards the interior of $S$:
\begin{equation}
H^i:=\{\mathbf{x} \in \mathbb{H}^n |  \mathbf{x}~\mbox{\boldmath$b$}^i =0 \} \ \ \text{with} \ \
\langle \mbox{\boldmath$b$}^i,\mbox{\boldmath$b$}^i \rangle = 1.
\end{equation}

The Grammian matrix $G(\mathcal{S}):=( \langle \mbox{\boldmath$b$}^i,
\mbox{\boldmath$b$}^j \rangle )_{i,j} ~ {i,j \in \{ 0,1,2 \dots n \} }$  is an
indecomposable symmetric matrix of signature $(1,n)$ with entries
$\langle \mbox{\boldmath$b$}^i,\mbox{\boldmath$b$}^i \rangle = 1$
and $\langle \mbox{\boldmath$b$}^i,\mbox{\boldmath$b$}^j \rangle
\leq 0$ for $i \ne j$ where

$$
\langle \mbox{\boldmath$b$}^i,\mbox{\boldmath$b$}^j \rangle =
\left\{
\begin{aligned}
&0 & &\text{if}~H^i \perp H^j,\\
&-\cos{\alpha^{ij}} & &\text{if}~H^i,H^j ~ \text{intersect \ along an edge of $P$ \ at \ angle} \ \alpha^{ij}, \\
&-1 & &\text{if}~\ H^i,H^j ~ \text{are parallel in the hyperbolic sense}, \\
&-\cosh{l^{ij}} & &\text{if}~H^i,H^j ~ \text{admit a common perpendicular of length} \ l^{ij}.
\end{aligned}
\right.
$$

For the {\it complete Coxeter orthoschemes} $\mathcal{S} \subset \overline{\mathbb{H}}^n$, 
we adopt the conventions for the Coxeter graph:
if the dihedral angle between two hyperplane faces is $\frac{\pi}{m_{ij}}$, then two related nodes are joined by an edge of weight $m_{ij}$. If the two hyperplane faces are orthogonal then the two corresponding nodes are disjoint.
In the hyperbolic case if two bounding hyperplanes of $S$ are parallel, 
then the corresponding nodes
are joined by a line marked $\infty$. 
If they are divergent then their nodes are joined by a dotted line. 

In general, the Coxeter orthoschemes were 
classified by {H-C. Im Hof} in \cite{IH1} and \cite{IH2} where he proved that 
they exist in dimension $\leq 9$, and gave a full list. 

{\it In this paper we consider a class of simply truncated orthoschemes where
$d=1$. That means, that one of its principal vertices (e.g. $A_{n-1}$) 
is the outer point related to the projective model, 
so they are truncated by its polar hyperplane
$Pol(A_{n-1})$. Moreover, we assume that $Pol(A_{n-1})$ and the hyperplane $H_0$ are parallel.
We obtain the above classification of the Coxeter orthoschemes that 
such tilings exist only in $3,4,5$ and $6$-dimensional hyperbolic spaces $\overline{\mathbb{H}}^n$.}

Packing problems related to the three-dimensional case were investigated 
in papers \cite{YSz22}, \cite{ASz} thus in 
the following we concentrate on $\overline{\mathbb{H}}^n$ with dimensions $n=4,5,6$.
We note here that with these analyses, we close the examination of these tiling types 
and the corresponding packing problems.

\begin{figure}[h!]
\centering
\includegraphics[width=15cm]{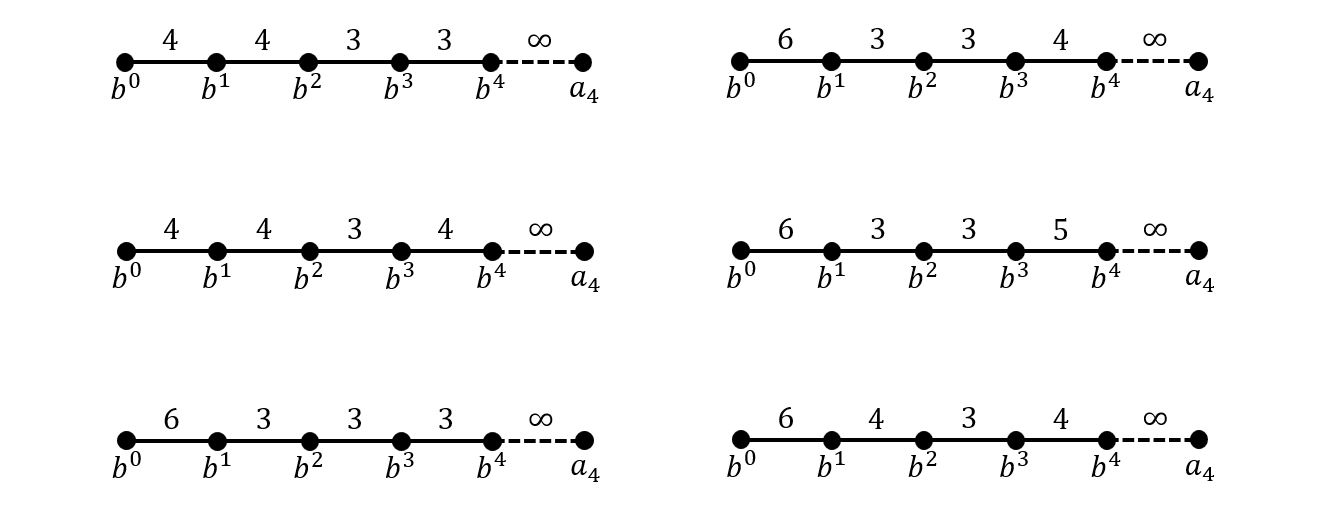}
\caption{Coxeter graph of simply truncated orthoschemes with parallel faces in $\overline{\mathbb{H}}^4$}
\label{Coxeter4}
\end{figure}

We begin by recalling the list of Coxeter graphs that give 
truncated orthoscheme with parallel faces (see  Fig.\ref{Coxeter4}--\ref{Coxeter56}),
and describing their Coxeter graph for $n=4,5,6$.
\begin{figure}[h!]
\centering 
\includegraphics[width=13cm]{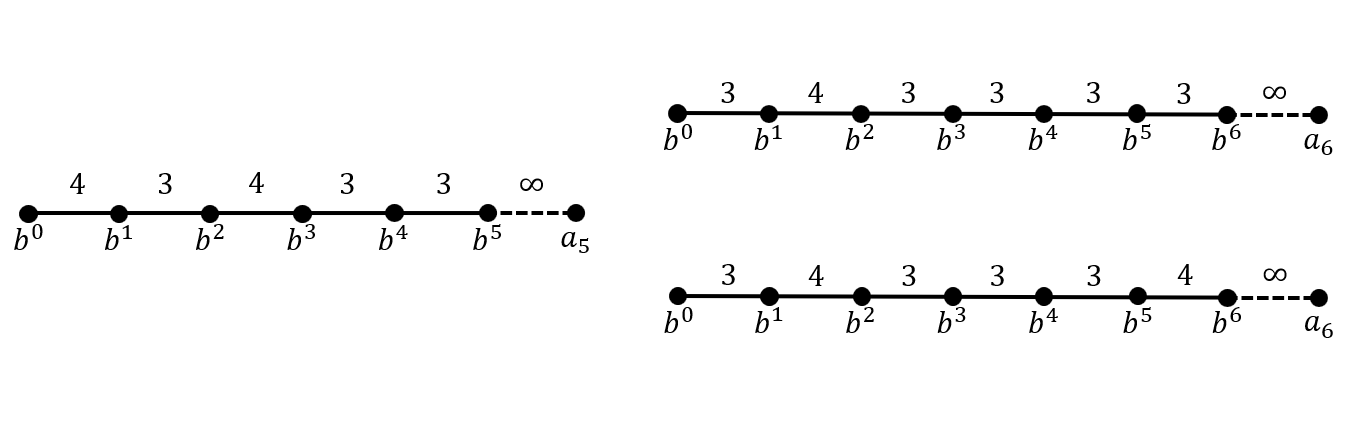} 
\caption{Coxeter graph of simply truncated orthoschemes with parallel faces 
in $\overline{\mathbb{H}}^5$ (left) and in $\overline{\mathbb{H}}^6$ (right)}
\label{Coxeter56} 
\end{figure}
We have obtained $6$ types of Coxeter tilings in $\overline{\mathbb{H}}^4$, $1$ type and $2$ types 
for $\overline{\mathbb{H}}^5$ and $\overline{\mathbb{H}}^6$ respectively. The structures of truncated 
orthoschemes (the fundamental domains of the considered tilings) 
in $\overline{\mathbb{H}}^4$ and $\overline{\mathbb{H}}^5$ are illustrated in Fig.~4.
\begin{figure}[h!]
\centering 
\includegraphics[width=12cm]{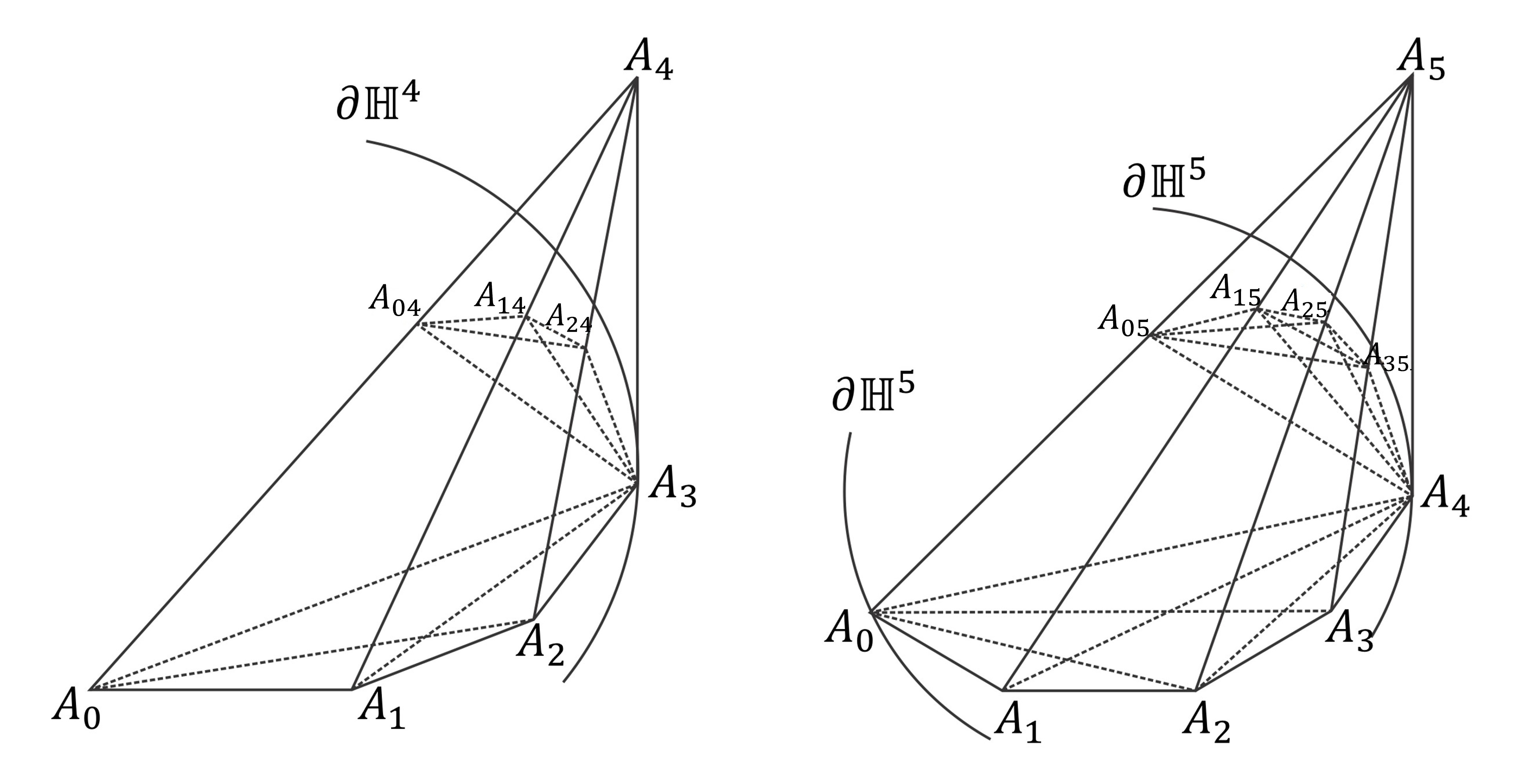}
\caption{The structure of a truncated orthoscheme in $\overline{\mathbb{H}}^4$, 
the vertex $A_3$ is an ideal, $A_4$ is ultra-ideal, while $A_0$ is either proper 
or ideal.
The polar hyperplane of $A_4$ is given by $A_3 A_{04} A_{14} A_{24}$
(left). The truncated orthoscheme in $\overline{\mathbb{H}}^5$, 
the vertex $A_4$ is an ideal, $A_5$ is ultra-ideal, 
while $A_0$ is either proper or ideal. 
The polar hyperplane of $A_5$ is represented by points $A_4 
A_{05} A_{15} A_{25} A_{35}$
(right)}
\label{Coxeter56} 
\end{figure}
Therefore, we shall consider the classical ball and horoball packings related to 
truncated orthoschemes with these Schl\"{a}fli symbols. 
We shall apply the methods used in $3$-dimensional cases (\cite{YSz22,ASz} with some 
important generalizations to higher dimensions.
\subsubsection{The volumes of orthoschemes in even dimensions}
The area formula for planar orthoschemes, which are right-angled triangles, is widely known as the defect formula. L. Schl\"afli extended this formula to spherical orthoschemes of even dimensions. Schl\"afli's reduction formula expresses the volume of an even-dimensional spherical orthoscheme in terms of the volumes of certain lower-dimensional orthoschemes. This formula can be straightforwardly extended to the hyperbolic case through analytic continuation.

In \cite{KR1}, R. Kellerhals further extended this formula to a specific class of hyperbolic polytopes called (complete) orthoschemes of degree $d$ (where $d \in {0, 1, 2}$). These polytopes arise in the study of hyperbolic Coxeter groups and are a particular type of fundamental polytope. R. Kellerhals demonstrated that the generalized reduction formula holds for even-dimensional complete orthoschemes and determined their volumes for all possible complete orthoschemes in dimensions $n$, where $n \in {4, 6, 8}$.

In Table \ref{Vol46}, we recall the volumes obtained for the even-dimensional 
orthoschemes we examined. 
\begin{table}[h!]
	\centering
		\begin{tabular}{||c c c||}
		\hline\xrowht[()]{10pt} 
		Dimension & Schl\"afli symbol & $\mathrm{vol}(\widehat{\mathcal{S}})$.  \\ 
		\hline\hline\xrowht[()]{10pt} 
		4&$\{4,4,3,3,\infty\}$ & $\frac{\pi^2}{288} \approx0.03427$  
		\\ 
		\hline\xrowht[()]{10pt}
		4&$\{4,4,3,4,\infty\}$ &  $ \frac{\pi^2}{144} \approx0.06854$ 
		\\
		\hline\xrowht[()]{10pt}
		4&$\{6,3,3,3,\infty\}$ &  $ \frac{\pi^2}{540} \approx 0.01828$ 
		\\
		\hline\xrowht[()]{10pt}
		4&$\{6,3,3,4,\infty\}$ & $ \frac{\pi^2}{288} \approx 0.03427$ 
		\\
		\hline\xrowht[()]{10pt}
		4&$\{6,3,3,5,\infty\}$ &  $ \frac{61\pi^2}{900} \approx 0.66894 $
		\\
		\hline\xrowht[()]{10pt}
		4&$\{6,3,4,3,\infty\}$ &  $ \frac{5\pi^2}{864} \approx 0.05712$ 
		\\
		\hline\hline\xrowht[()]{10pt}
		6&$\{3,4,3,3,3,3,\infty\}$ & $\frac{11 \pi^3}{86,400} \approx0.00395$  
		\\
		\hline\xrowht[()]{10pt}
		6&$\{3,4,3,3,3,4,\infty\}$ &  $ \frac{ \pi^3}{86,400} \approx 0.00036$ 
		\\
		\hline
	\end{tabular}
\caption{Volumes of truncated orthoschemes with parallel faces in $\overline{\mathbb{H}}^4$ and $\overline{\mathbb{H}}^6$.}
\label{Vol46}
\end{table}
\subsubsection{Volumes of orthoschemes in $\overline{\mathbb{H}}^5$}
In contrast to the cases of $\overline{\mathbb{H}}^4$ and $\overline{\mathbb{H}}^6$, the volumes of simply truncated orthoschemes in $\overline{\mathbb{H}}^5$ cannot be expressed by a simple closed formula. To overcome this challenge, we employ a method known as the scissor congruence method, which involves decomposing the orthoscheme into several smaller orthoschemes (see \cite{KR2}). These smaller orthoschemes must satisfy the condition of being double asymptotic orthoschemes, allowing their volumes to be computed using a formula established by R. Kellerhals in \cite{KR3, KR4}. This formula utilizes the ellipticity and parabolicity conditions of the principal vertices, specifically $A_0$ and $A_5$, within the orthoscheme $\mathcal{S}$. \\
Consider an orthoscheme $\mathcal{S}$ with a Coxeter-Schl\"afli matrix denoted as $(b^{ij})$, and let ${v_1, v_2, v_3, v_4, v_5}$ represent its Schl\"{a}fli symbol. Within this context, $(b^{ij})_0$ and $(b^{ij})_5$ refer to the principal submatrices of $(b^{ij})$ obtained by removing the corresponding rows and columns associated with the vertices $A_0$ and $A_5$ respectively.
Geometrically, this removal process corresponds to eliminating the opposite face of the vertex $A_0$ or $A_5$ from the orthoscheme $\mathcal{S}$, resulting in a simplicial cone structure.

\begin{figure}[h!]
\centering 
\includegraphics[scale=0.40]{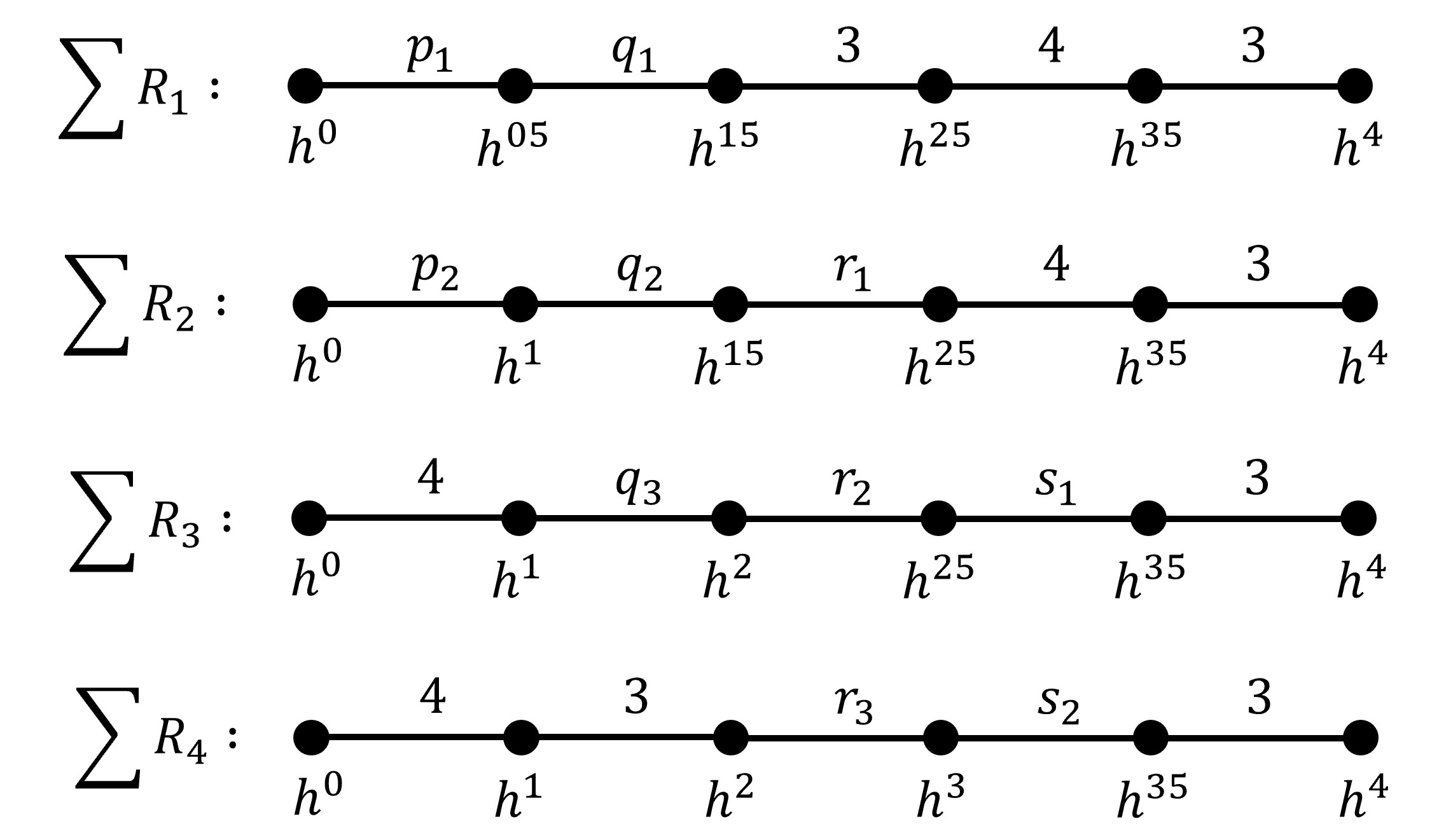}  
\caption{Coxeter graph of dissected doubly asymptotic orthoschemes $R_j$, $(j \in \{1,2,3,4\})$. 
Here, $h^{k}$ denotes the opposite face of vertex $A_k$ in each corresponding $R_j$}
\label{Dissected5}
\end{figure}

We have the following elliptic/parabolic conditions:
\begin{enumerate}
    \item If $(b^{ij})_0$ and $(b^{ij})_5$ are positive definite, then $A_0$ and $A_5$ are a proper vertices and it said to be \textit{elliptic}.
    \item If $(b^{ij})_0$ and $(b^{ij})_5$ are positive semi-definite, then $A_0$ and $A_5$ are ideal vertices and it said to be \textit{parabolic}.
\end{enumerate}
Through the analysis of the corresponding polynomials associated with these submatrices, the following findings emerge:
\begin{enumerate}
    \item $(b^{ij})_0$, and $(b^{ij})_5$ are positive definite, \\ $\displaystyle \frac{\cos^2{\left( \frac{\pi}{v_3} \right)}}
    {\sin^2{\left( \frac{\pi}{v_2} \right)}}+\frac{\cos^2{\left( \frac{\pi}{v_4} \right)}}{\sin^2{\left( \frac{\pi}{v_5} \right)}}<1$, \ \ 
 $\displaystyle \frac{\cos^2{\left( \frac{\pi}{v_2} \right)}}{\sin^2{\left( \frac{\pi}{v_1} \right)}}+\frac{\cos^2{\left( \frac{\pi}{v_3} 
 \right)}}{\sin^2{\left( \frac{\pi}{v_4} \right)}}<1$.
  \item $(b^{ij})_0$ and $(b^{ij})_5$ are semi-positive definite, \\ $\displaystyle \frac{\cos^2{\left( \frac{\pi}{v_3} \right)}}
  {\sin^2{\left( \frac{\pi}{v_2} \right)}}+\frac{\cos^2{\left( \frac{\pi}{v_4} \right)}}{\sin^2{\left( \frac{\pi}{v_5} \right)}}=1$, \ \    
 $\displaystyle \frac{\cos^2{\left( \frac{\pi}{v_2} \right)}}{\sin^2{\left( \frac{\pi}{v_1} \right)}}+\frac{\cos^2{\left( \frac{\pi}{v_3} 
 \right)}}{\sin^2{\left( \frac{\pi}{v_4} \right)}}=1$.
\end{enumerate}
By employing the scissor congruence method, we can decompose the simply truncated $\widehat{\mathcal{S}}$ orthoscheme with the Schl\"afli symbol $\{4,3,4,3,3,\infty\}$ into several smaller orthoschemes.\\
The polar hyperplane $Pol(A_5)$ of the ultra-ideal vertex $A_5$ encompasses the ideal vertex $A_4$ of $\mathcal{S}$. The intersections of $Pol(A_5)$ with the edges $A_iA_5$ (where $i\in\{0,1,2,3\}$) are denoted as $A_{i5}$. For a visual representation, please refer to Figure 3.
The truncated orthoscheme $\widehat{\mathcal{S}}$ can be dissected into 
some ``smaller" doubly asymptotic orthoschemes $R_j$, $(j \in \{1,\dots 4\})$ 
where each of them is the convex hull of vertices:
\begin{align*}
    &R_1:=conv(A_0,A_{05},A_{15},A_{25},A_{35},A_4)~~~&R_2:=conv(A_0,A_{1},A_{15},A_{25},A_{35},A_4) \\
    &R_3:=conv(A_0,A_{1},A_{2},A_{25},A_{35},A_4)~~~~&R_4:=conv(A_0,A_{1},A_{2},A_{3},A_{35},A_4).
\end{align*}
Therefore, the volume of $\widehat{\mathcal{S}}$ can be computed as
\begin{equation*}
    \mathrm{vol}(\widehat{\mathcal{S}})=\sum_{i=1}^{4} \mathrm{vol}(R_i).
\end{equation*}
By following the method described by Kellerhals in \cite{KR2}, 
the dihedral angles of orthoschemes $R_1, R_2, R_3$, and $R_4$ 
can be computed by the following method using their Coxeter graphs, (see Fig.~\ref{Dissected5}). 

Based on the structure of $\widehat{\mathcal{S}}$, we can derive the following two equations directly:
\begin{align*}
    \frac{\pi}{p_1}+\frac{\pi}{p_2}=\frac{\pi}{2}, ~~\frac{\pi}{s_1}+\frac{\pi}{s_2}=\frac{\pi}{2}.
\end{align*}
The first equation arises from the consideration of the angles $\frac{\pi}{p_1}$ and $\frac{\pi}{p_2}$ associated with the vertices $A_{15}, A_{25}, A_{35}, A_4$ in 3-space. These angles are present in the shared portion of the hyperplane faces $A_{05}A_{15}A_{25}A_{35}A_4$, $A_{0}A_{15}A_{25}A_{35}A_4$, and $A_{1}A_{15}A_{25}A_{35}A_4$.
The second equation results from the presence of the angles $\frac{\pi}{s_1}$ and $\frac{\pi}{s_2}$ on the edge $A_2A_4$, which lies in the 3-spaces determined by the vertices $A_2, A_3, A_{35}, A_4$ and $A_2A_{25}A_{35}A_4$.
By incorporating the parabolicity conditions and the above equations related to the parameters $p_1, p_2, s_1, s_2$, we can determine all the parameters in the Schl\"afli symbols depicted in Fig.~\ref{Dissected5}.
\begin{align*}
    p_1=p_2=s_1=s_2=4, ~~~~~~q_1=q_2=q_3=r_1=r_2=r_3=3.
\end{align*}
Finally, we determine that the Schl\"afli symbols for all the doubly asymptotic orthoschemes $R_j$ are ${4,3,3,4,3}$, and their volumes can be computed using the following formulas provided by R. Kellerhals.
\begin{prop}[\cite{KR3, KR4}]
Denote by $R \subset \overline{\mathbb{H}^5}$ a doubly asymptotic $5$-orthoscheme
represented by dihedral angles $\alpha_1$, $\alpha_2$, $\alpha_3$, $\alpha_4$, $\alpha_5$ with $\lambda=\tan{\theta}=\frac{\sqrt{|\mathrm{det}{(b^{ij})}|}}{\cos{\alpha_1}\cos{\alpha_3}\cos{\alpha_5}}$, $0 \leq \theta \leq \frac{\pi}{2}$. Let $0 \leq \alpha_0 \leq \frac{\pi}{2}$ such that $\tan{\alpha_0}=\cot{\theta} \tan{\alpha_3}$. Then,
\begin{align*}
    \mathrm{vol}_5(R)&=-\frac{1}{8}\{ I(\lambda^{-1},0;\alpha_1) + \frac{1}{2}I(\lambda,0;\alpha_2)-I(\lambda^{-1},0;\alpha_0')+
    \frac{1}{2}I(\lambda,0;\alpha_4)+\\ & +I(\lambda^{-1},0;\alpha_5) \}
    + \frac{1}{32}\{ I_{alt} (\lambda, \alpha_1; \alpha_2) + I_{alt} (\lambda, \alpha_5; \alpha_2) \}
\end{align*}
where $I(a,b;x)$ is the trilogarithmic function, i.e 
\begin{align*}
    I(a,b;x)=\int_{\pi/2}^{x} \frac{1}{2} \mathrm{Im} \left( \mathrm{Li}_2 (e^{2iy})  \right) d \left( a~\tan^{-1}(b+y) \right)
\end{align*}
with
\begin{align*}
    \mathrm{Li}_2(z)=-\int_{0}^{z} \frac{\mathrm{ln}(1-t)}{t} dt.
\end{align*}
and
\begin{align*}
    I_{alt} (a,b;x)&= I_{\delta}(a,b;x) + I_{\delta}(a,-b;x)\\
    I_{\delta}(a,b;x) &= I(a,-(\frac{\pi}{2}+b);\frac{\pi}{2}+b+x)-I(a,-(\frac{\pi}{2}+b);\frac{\pi}{2}+b+\frac{\pi}{2})\\
    \alpha' &=\frac{\pi}{2}-\alpha.
\end{align*}
\end{prop}
In the present cases, each doubly asymptotic orthoscheme $R_j$ possesses the Schl\"afli symbol $\{4,3,3,4,3\}$, which implies $\lambda=1$ and $\cos^2{\alpha_1}+\cos^2{\alpha_2}+\cos^2{\alpha_3}=1$. Therefore, we can apply the specific version of the proposition as follows
\begin{prop}[\cite{KR3, KR4}]\label{Vol5}
    Let $R$ denote a doubly asymptotic 5-orthoscheme with $\lambda=\tan{\theta}=1$ that is
    \begin{equation*}
        \{\alpha_1, \alpha_2, \alpha_3, \alpha_4, \alpha_5 \}, ~\text{with}~\cos^2{\alpha_1}+\cos^2{\alpha_2}+\cos^2{\alpha_3}=1
    \end{equation*}
    Then,
    \begin{align*}
        \mathrm{vol}(R)&=\frac{1}{4}\biggl\{\mathtt{J}_3(\alpha_1)+\mathtt{J}_3(\alpha_2)-\frac{1}{2}\mathtt{J}_3\left(\frac{\pi}{2}-\alpha_3\right)\biggr\}\\
        &-\frac{1}{16}\biggl\{ \mathtt{J}_3 \left(\frac{\pi}{2}+\alpha_1+\alpha_2\right)+ \mathtt{J}_3 \left(\frac{\pi}{2}-\alpha_1+\alpha_2\right) \biggr\}+\frac{3}{64}\zeta{(3)},\\
        &\text{where},\\
        &\mathtt{J}_3(\alpha)=\frac{1}{4}\mathrm{Re}\left( \mathrm{Li}_3\left( e^{2i\alpha} \right) \right),~\text{(the Lobachevsky function of order 3)}\\
        &\text{with},\mathrm{Li}_k(z)=\int_{0}^{z}\frac{\mathrm{Li}_{k-1}(t)}{t}dt~~\text{and}~~ \mathrm{Li}_1(z)=-\mathrm{log}(1-z),\\
        &\text{and}~ \zeta~\text{is Riemann Zeta function}.
    \end{align*}
\end{prop}
The volume of each doubly asymptotic orthoscheme $R_j$, where $j \in {1,\dots,4}$, can be computed using the formula provided in Proposition \ref{Vol5}. It is known that this volume is equal to $\frac{7\zeta(3)}{4608}$, as stated in \cite{KR3}. Therefore, the volume of $\widehat{\mathcal{S}}$ can be evaluated directly as follows:
\begin{align*}
    \mathrm{vol}(\widehat{\mathcal{S}})&=\sum_{i=1}^{4} \mathrm{vol}(R_i)=4\cdot \frac{7\zeta{(3)}}{4608}=\frac{7\zeta{(3)}}{1152}.
\end{align*}
\section{Packing with classical balls in $\overline{\mathbb{H}}^n$} \label{ClassicBall3}
To construct and determine the inballs and their radii, as well as the optimal ball packing density related to the investigated Coxeter tilings, we can directly apply the method described in \cite{J14, YSz22}, along with the formulas (\ref{inradiusConditions}), (\ref{Inradius}), and (\ref{NS-condition}). When determining the inradius of truncated orthoschemes, we need to consider two distinct types. The first type involves situations where the inscribed ball of a complete orthoscheme is the same as the inscribed ball of the truncated orthoscheme. The second type pertains to situations where this equality does not hold.
\begin{enumerate}
\item {\bf Type 1}

M.~Jacquemet determined the inradii of truncated simplices in $n$-dimensional hyperbolic spaces when their inballs do not have common inner points with the corresponding truncating hyperplanes (see \cite{J14}). Now, let us recall some important statements from this remarkable paper.
We will use the denotations introduced in the previous section. Let $\mathcal{S}$ be a complete orthoscheme and $\widehat{\mathcal{S}}$ be the truncated orthoscheme under consideration. The corresponding Coxeter-Schl\"afli matrix is denoted as $(B^{ij})$, with its principal submatrix denoted as $(b^{ij})$, which is the Coxeter-Schl\"afli matrix of $S$. The following lemma provides a sufficient condition for the existence of an inball in a truncated simplex.
\begin{lem}[\cite{J14}]\label{inradiusConditions}
A truncated hyperbolic simplex $\widehat{\mathcal{S}}$ with Coxeter-Schl\"afli principal submatrix $(b^{ij})$ has inball, 
(embedded ball of maximal finite radius) in $\mathbb{H}^n$ if and only 
if $\sum_{i,j=1}^{n+1} \mathrm{cof}_{ij}({b^{ij}}) > 0$.
\end{lem}
It can be reformulated using the matrix $(a_{ij})= (b^{ij})^{-1}$ by the following
\begin{equation}
\sum_{i,j=1}^{n+1} \mathrm{det}(b^{ij})(a_{ij}) >0.
\end{equation}
Moreover, the following lemma states the formula of inradius in a complete orthoscheme.
\begin{lem}[\cite{J14}]\label{Inradius}
Let $(b^{ij})$ be the Coxeter--Schl\"afli matrix of complete orthoscheme ${\mathcal{S}}$ with inball ${\mathcal{B}} \subset \overline{\mathbb{H}}^n$. 
Then, the inradius $\mathtt{r}=\mathtt{r}({\mathcal{B}})$ is given by
\begin{equation}
\mathtt{r}=\sinh^{-1}{\sqrt{-\frac{1}{\sum_{i,j=1}^{n+1} a_{ij}}}}.
\end{equation}
\end{lem}
In our investigation, the orthoscheme contains one ultra-ideal vertex located outside the Beltrami-Cayley-Klein model. Consequently, the orthoscheme is truncated by a polar hyperplane associated with its ultra-ideal vertex. Thus, it is important to examine whether the inradius of the original orthoscheme and the inradius of the truncated orthoscheme are equal. Lemma \ref{NS-condition} provides the necessary and sufficient conditions for the coincidence of the inradius between the truncated simplex (orthoscheme) and the complete simplex (orthoscheme).
\begin{lem}[\cite{J14}]\label{NS-condition}
Let ${\mathcal{S}}$ be an $n$-dimensional simplex with $k \le n$ vertices are ultra ideal, with Coxeter-Schl\"afli matrix $(b^{ij})$, such that ${\mathcal{S}}$ has an inball 
$\mathcal{B} \subset \mathbb{H}^n$ of radius $\mathtt{r}$. Denote by $\widehat{\mathcal{S}} \subset \overline{\mathbb{H}^n}$ its associated hyperbolic $k$-truncated simplex 
with respect to the ultra-ideal vertices $v_1, \cdots, v_k$, $1 \leq k$. Let $\widehat{\mathtt{r}}$ be the inradius of inball of $\widehat{\mathcal{S}}$. Then, $\mathtt{r}=\widehat{\mathtt{r}}$ if and only if
\begin{equation}
    \frac{\sum_{j=1}^{n+1}\mathrm{cof}_{ij}{(b^{ij})}}{\mathrm{det}({b^{ij}})\mathrm{cof}_{ii}{(b^{ij})}} \geq 1 ~~~ \mathrm{for~all}~ i=1, \cdots, k.
\end{equation}
\end{lem}
\item {\bf Type 2}

We consider the case where the inradius of the complete orthoschemes and the truncated orthoschemes is not the same, i.e., when Lemma \ref{NS-condition} does not hold. In these cases, the constructed inball intersects the truncation hyperplane:
\begin{align*}
\mathtt{r}=d(\bc,H^0)=\cdots=d(\bc,H^{n1}) > d(\bc, H^{n+1}),
\end{align*}
where $\bc$ is the insphere center of the complete orthoscheme $\mathcal{S}$. We applied the classical way to determine the incenter and the radius in the mentioned cases. 
In general, to find the center and the radius $\widehat{\mathtt{r}}$ of optimal inball we determine the hyperplane bisectors of faces of truncated orthoscheme. 
The inball of the maximal radius has to touch at least four faces of $\widehat{\mathcal{S}}$ and one of them must be the face determined by form $[\Bb^n]$. 
Therefore, we have $n+1$ analogs cases that provide the candidates for the optimal incenter. We have to determine these centers and 
select the center with the maximum radius. This optimal ball is denoted by $\mathcal{B}$. 
\end{enumerate}
We introduce the local density function $\delta_{opt}(\widehat{\mathcal{S}})$ related to orthoscheme $\widehat{\mathcal{S}}$ generated tiling:
\begin{defn}
The local density function $\delta_{opt}(\widehat{\mathcal{S}})$ related to orthoscheme tilings generated by truncated orthoscheme $\widehat{\mathcal{S}}$:
\begin{equation}
\delta_{opt}(\widehat{\mathcal{S}}):=\frac{\mathrm{vol}(\mathcal{B})}{\mathrm{vol}(\widehat{\mathcal{S}})}. \notag
\end{equation}
\end{defn}
We obtain the volumes of the $\mathrm{vol}(\mathcal{B})$ balls by the classical formula 
$$\mathrm{vol}(\mathcal{B})=k^n \int_\Phi \sinh^{n-1}\phi^1\sin^{n-2}\phi^2 \dots \sin\phi^{n-1} \mathrm{d} \phi^1 \wedge \dots \mathrm{d} \phi^n, $$ 
where at present $k=1$ and
$$
\Phi: 0\le \phi^1 \le \frac{\mathtt{r}}{k},0\le \phi^2 <\pi, \dots, 0\le \phi^n <2\pi.
$$
The volumes of truncated simplices $\widehat{\mathcal{S}}$ can be calculated by the results of subsections 2.3.1-2. With the help of the procedure described above, we determined the densest 
ball packing configurations and their densities belonging to the investigated Coxeter tilings (the calculations are not detailed here, see \cite{YSz22,ASz} and \cite{J14}).
The results are summarized in the following tables:
\begin{table}[h!]
	\centering
\begin{footnotesize}
	\begin{tabular}{||c c c c c||} 
		\hline\xrowht[()]{10pt}
		Schl\"afli symbol & Inradius $\widehat{\mathtt{r}}$ & $\mathrm{vol}(\mathcal{B})$ & $\mathrm{vol}(\widehat{\mathcal{S}})$ & $\delta_{opt}(\widehat{\mathcal{S}})$ \\ [0.5ex] 
		\hline\hline\xrowht[()]{10pt}
		$\{4,4,3,3,\infty \}$ & $\approx 0.14440$  & $\approx 0.00216$ & $\frac{\pi^2}{288} \approx 0.03427$
		& $\approx \bold{0.06304}$
		\\ 
		\hline\xrowht[()]{10pt}
		$\{4,4,3,4,\infty \}$ & $\approx 0.15986$ & $\approx 0.00325$ & $\frac{\pi^2}{288} \approx 0.06854$
		& $\approx 0.04742$
		\\
		\hline\xrowht[()]{10pt}
		$\{6,3,3,3,\infty \}$ & $\approx 0.11182$ & $\approx 0.00077$
		& $\frac{\pi^2}{540} \approx 0.01828$
		& $\approx 0.04239$ \\
		\hline\xrowht[()]{10pt}
		$\{6,3,3,4,\infty \}$ & $\approx 0.12751$ & $\approx 0.00131$
		& $\frac{\pi^2}{288} \approx 0.03427$
		& $\approx 0.03828$ \\
		\hline\xrowht[()]{10pt}
		$\{6,3,3,5,\infty \}$ & $\approx 0.14200$ & $\approx 0.00202$
		& $\frac{61\pi^2}{900} \approx 0.66894$
		& $\approx 0.00302$
		\\
		\hline\xrowht[()]{10pt}
		$\{6,3,4,3,\infty \}$ & $\approx 0.16208$ & $\approx 0.00344$
		& $\frac{5\pi^2}{864} \approx 0.05712$ & $\approx 0.06014$
		\\ 
		\hline
		
	\end{tabular}
\end{footnotesize}
\caption{Optimal ball packings and their densities in $\overline{\mathbb{H}}^4$} 
\end{table}
 \begin{thm}
 In hyperbolic space $\overline{\mathbb{H}}^4$, between congruent ball packings of {\it classical balls}, generated by simply truncated Coxeter orthoschemes with parallel faces, 
 the $\mathcal{B}_{\{4,4,3,3\}}$ ball configuration provides the densest packing with density 
 $\approx 0.06304 $. 
 \end{thm}
 \begin{table}[h!]
 	\centering
 \begin{footnotesize}
 	\begin{tabular}{||c c c c c||} 
 		\hline\xrowht[()]{10pt}
 		Schl\"afli symbol & Inradius $\widehat{\mathtt{r}}$ & $\mathrm{vol}(\mathcal{B})$ & $\mathrm{vol}(\widehat{\mathcal{S}})$ & $\delta_{opt}(\widehat{\mathcal{S}})$ \\ [0.5ex] 
 		\hline\hline\xrowht[()]{10pt}
 		$\{4,3,4,3,3,\infty \}$ & $\approx 0.11414$  & $\approx 0.00010$
 		& $\frac{5}{7\zeta{(3)}}$
 		& $\approx \bold{0.0140}$
 		\\ 
 		\hline
 		\end{tabular}
 \end{footnotesize}
 \caption{Optimal ball packings and their densities in $\overline{\mathbb{H}}^5$} 
\end{table}
 \begin{thm}
 In hyperbolic space $\overline{\mathbb{H}}^5$, between congruent ball packings of {\it classical balls}, generated by simply truncated Coxeter orthoschemes with parallel faces, 
 the $\mathcal{B}_{\{4,3,4,3,3\}}$ ball configuration provides the densest packing with density $\approx 0.0140 $. 
 \end{thm}
 \begin{table}
\begin{footnotesize}
	\begin{tabular}{||c c c c c||} 
		\hline\xrowht[()]{10pt}
		Schl\"afli symbol & Inradius $\widehat{\mathtt{r}}$ & $\mathrm{vol}(\mathcal{B})$ & $\mathrm{vol}(\widehat{\mathcal{S}})$ & $\delta_{opt}(\widehat{\mathcal{S}})$ \\ [0.5ex] 
		\hline\hline\xrowht[()]{10pt}
		$\{3,4,3,3,3,3,\infty \}$  & $\approx 0.06102$
		& $\approx 2.67258 \times 10^{-7}$ & $\frac{11 \pi^3}{86,400} \approx 0.00395$ & $\approx 6.77032 \times 10^{-5}$
		\\ 
		\hline\xrowht[()]{10pt}
		$\{3,4,3,3,3,4,\infty\}$ & $\approx 0.06102$ & $\approx 2.67258 \times10^{-7}$ & $\frac{ \pi^3}{86,400} \approx 0.00036$
		& $ \approx \bold{7.44867} \times \bold{10^{-4}}$
		\\
		\hline
	\end{tabular}
\caption{Optimal ball packings and their densities in $\overline{\mathbb{H}}^6$}
\end{footnotesize}
\end{table}
 \begin{thm}
 In hyperbolic space $\overline{\mathbb{H}}^6$, between congruent ball packings of {\it classical balls}, generated by simply truncated Coxeter orthoschemes with parallel faces, 
 the $\mathcal{B}_{\{3,4,3,3,3,4,\infty\}}$ ball configuration provides the densest packing with density $\approx 7.44867 \times 10^{-4}$. 
 \end{thm}
\section{Horoball packings in $\overline{\mathbb{H}}^n$, $n=4,5,6$}
\subsection{Horoball packing density and important lemmas}
In this subsection, we define packing density and collect three 
Lemmas used in the next section to find the optimal packing 
densities for the examined orthoscheme tilings. 

Let $\cT$ be a Coxeter orthoscheme tiling of $\overline{\mathbb{H}}^n$ given by Schl\"afli symbols in Fig.~1-2.
The symmetry group of a Coxeter tiling $\Gamma_\cT$ 
contains its Coxeter group and isometric mapping between two cells 
in $\cT$ preserves the tiling.
Any simplex cell of $\cT$ acts as a fundamental domain $\cF_{\cT}$
of $\Gamma_\cT$, and the Coxeter group is generated by reflections 
on the $(n - 1)$-dimensional facets of $\cF_{\cT}$. 
In this paper, we consider only noncompact simply truncated Coxeter orthoscemes 
(and the corresponding tilings) 
with parallel faces 
with one or more ideal vertex, 
then the orbifold $\overline{\mathbb{H}}^n/\Gamma_{\cT}$ has at least one cusp (see Fig.~1-2)
for $4 \le n \le 6$. 

Define the density of a regular horoball packing 
$\mathcal{B}_{\cT}$ of Coxeter orthoscheme tiling $\cT$ as
\begin{equation}
\delta(\mathcal{B}_{\cT})=
\frac{\sum_{i=1}^m \mathrm{vol}(B_i \cap \cF_{\cT})}{\mathrm{vol}(\cF_{\cT})}.
\label{eq:density}
\end{equation}
$\cF_{\cT}$ denotes the fundamental domain of tiling $\cT$, 
$m$ the number of ideal vertices of $\cF_\cT$, 
and $B_i$ the horoball centered at the $i$-th ideal vertex. 
We allow horoballs of different types at each asymptotic vertex of the tiling. 
A particular set of horoballs $\{B_i\}_{i=1}^m$ with different horoball types is 
allowed if it gives a packing: no two horoballs may have a common interior point, 
and we require that no horoball extend beyond the facet opposite to the vertex where 
it is centered. The second condition ensures that the packing remains invariant under the actions of $\Gamma_\cT$ with $\cF_\cT$.
With these conditions satisfied, the packing density in $\cF_{\cT}$ extends 
to the entire $\overline{\mathbb{H}}^n$ by actions of $\Gamma_{\tau}$.
In the case of Coxeter truncated orthoscheme tilings, 
Dirichlet--Voronoi cells coincide with the fundamental domains (truncated orthoschemes).
We denote the optimal horoball packing density as
\begin{equation}
\delta_{opt}(\cT_\Gamma) = \sup\limits_{\mathcal{B}_{\cT} \text{~packing}} 
\delta(\mathcal{B}_{\tau}).
\label{eq:opt_density}
\end{equation}

Let $\{A_i\}_{i=0}^n \in \mathbb{P}^n$ denote the vertices of the starting orthoschemes of the 
considered Coxeter tiling $\cT_{\Gamma}$ where $A_n$ is an outer point (see Fig.~3) 
and the other vertices are proper or ideal points in the model. 
We obtain from this consideration the fundamental domain $\cF_{\Gamma}$ of $\cT_{\Gamma}$ by truncating the orthoscheme with 
the polar hyperplane $Pol(A_n)$ of $A_n$. In these considered cases (see Fig.~1-2), $\cF_{\Gamma}$
is an $n$-polytope with vertices $A_0A_1A_2\dots A_{n-1}A_{0n}A_{1n}\dots A_{(n-2)n}$ where
$A_{n-1}$ lies on the boundary $\partial \mathbb{H}^n$, $A_0$ is either proper or ideal vertex and the other vertices are proper points of 
$\mathbb{H}^n$ (see Fig.~3). In the case of the vertices lying at infinity, we particularly have two cases for the hyperball packings related to tilings.
\begin{enumerate}
\item {\bf One horoball types}
Here we describe a procedure for finding the optimal horoball packing density in the fundamental domain $\cF_{\Gamma}$ with a single ideal vertex $A_{n-1}$. 
Packing density is maximized by the largest horoball type admissible in cell $\cF_{\Gamma}$ centered at $A_{n-1}$. 

We transform $\cF_{\Gamma}$ using hyperbolic isometry that takes the vertex $A_{n-1}$ to the point $(1, 0, \dots, 0, 1)$.
The coordinates of the other vertices of $\cF_{\Gamma}$ are determined by the dihedral angles of 
$\cF_{\Gamma}$ indicated in the Coxeter diagrams in Fig.~1-2. 

Let $\mathcal{B}_{n-1}(s)$ denote the 1-parameter family of horoballs centered at $A_{n-1}$ where 
$s$-parameter related to the Busemann function measures the ``radius" of the horoball, the minimal Euclidean signed distance between the horoball 
and the center of the model $O$, it is taken to be negative if the horoball contains the model center.

The maximal horoball $\mathcal{B}_{n-1}(s)$ has to be tangent to
at least one of the hyperplanes of the fundamental domain $\cF_{\Gamma}$ that does not contain the vertex $A_{n-1}$, so that it does not intersect another such hyperplane.

The possible tangent points $[\mathbf{f}_i]$ of $\mathcal{B}_{n-1}(s)$ and the corresponding hyperplanes are determined by the projection of vertex $A_{n-1}$ on the possible hyperplanes 
given by forms $\Bu_i$ (see Fig.~4),
\begin{equation}
\mathbf{f}_i^{n-1} =\ba_{n-1} - \frac{\langle \ba_{n-1}, \mathbf{u}_i \rangle}{\langle \Bu_i,\Bu_i \rangle} 
\mathbf{u}_i.
\label{eq:u4fp}
\end{equation}
\begin{figure}[h!]
\centering 
\includegraphics[width=5cm]{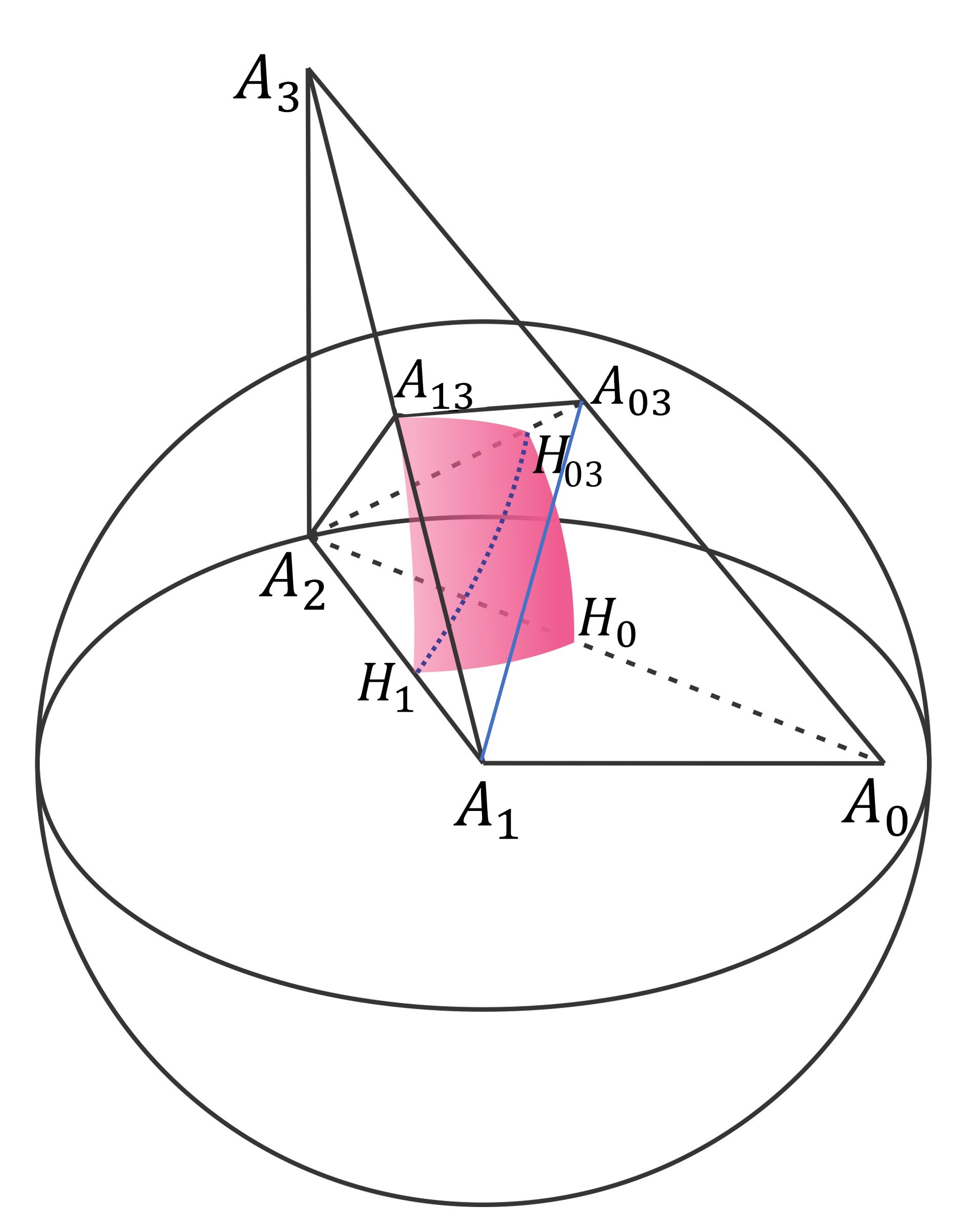} \includegraphics[width=6cm]{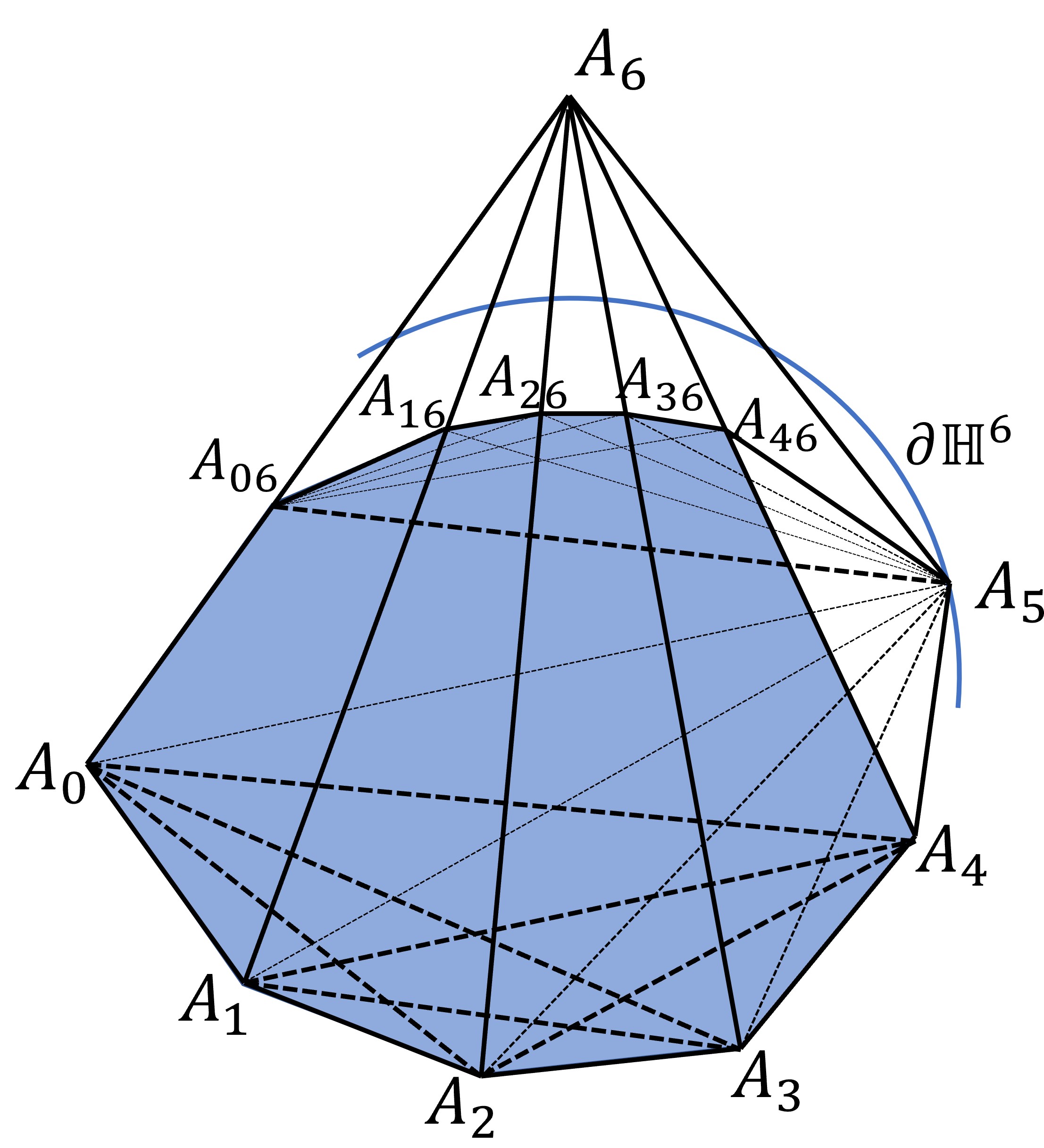} 
\caption{The horosphere $\mathcal{B}_{2}$ touch its opposite face and the obtained horospherical polyhedron ${A}_{13} {H}_0 {H}_1 {H}_{03}$ is dissected into two simplices 
${H}_0 {H}_1 {H}_{03}$ and ${A}_{13} {H}_1 {H}_{03}$ 
through hyperplane ${A}_2 {A}_{1} {A}_{03}$ in the $3$-dimensional hyperbolic space $\overline{\mathbb{H}}^3$ (left). The ``Opposite" faces of $A_5$ 
(in the shaded region), indicated by vertices 
$A_0$, $A_1$, $A_2$, $A_3$, $A_4$ $A_{06}$, $A_{16}$, $A_{26}$, $A_{36}$, $A_{46}$ in a truncated orthoscheme of $\overline{\mathbb{H}}^6$ (right).}
\label{HoroPolihedron3}
\end{figure}
The optimal value of the $s$-parameter for the maximal horoball can 
be determined from the above equations of the horosphere through 
$A_{n-1}$ and footpoints $\mathbf{f}_i^{n-1}$. 

The intersections $[\bh_i]$ of horosphere $\partial \mathcal{B}_{n-1}$ and the edges of the polytope $\cF_{\Gamma}$ are found by 
parameterizing the edges $\bh_i(\lambda) = \lambda \ba_{n-1} +\ba_i$ $(i \in \{ 0, \cdots, n-2, 0n, \cdots, (n-2)n \})$ then finding their intersections with $\partial \mathcal{B}_{n-1}$. 
The volume of the horospherical $(n-1)$-simplex determines the volume of the horoball piece by equation \eqref{eq:bolyai}.
The data for the horospherical $(n-1)$-polytope is obtained by finding hyperbolic distances $l_{ij}$ via equation \eqref{prop_dist},
$l_{ij} = d(H_i, H_j)$ where $d(\bh_i,\bh_j)= \arccos\left(\frac{-\langle \bh_i, \bh_j 
\rangle}{\sqrt{\langle \bh_i, \bh_i \rangle \langle \bh_j, \bh_j \rangle}}\right)$.
Moreover, the horospherical distances $L_{ij}$ are found by formula \eqref{eq:horo_dist}.

We obtain a horospherical $(n-1)$-polytope $\mathcal{P}$ whose edges are known.
The intrinsic geometry of a horosphere is Euclidean, 
so the obtained polytope's volume is computed in 
Euclidean sense. The computation of the volume of the above $n$-dimensional polytope is determined using its decomposition into $n$-dimensional 
simplexes. 
In $\overline{\mathbb{H}}^n$, the opposite faces of ideal vertex $A_{n-1}$ is determined by polytope $\mathcal{P}$ that is formed as a convex hull of 
$A_0,\cdots,A_{n-1},A_{0n},\cdots,A_{(n-2)n}$, see shaded region in Fig.~\ref{HoroPolihedron3}. 
This polytope $\mathcal{P}$ can be dissected into $n-1$ simplices as follows:
\begin{align*}
    E_k:=conv(A_{0n},\cdots A_{kn}, A_k,\cdots, A_{n-2}), ~~~\text{for } k=0,\cdots,n-2
\end{align*}
Note that,
\begin{align*}
    {int}(E_{k_1}) \cap {int}(E_{k_2})&=\emptyset, \text{for}~k_1\neq k_2,~\text{where} ~int(E_k)~ \text{denote the interior of} E_k\\
    \text{and} \bigcup_{k=0}^{n-2}E_k&=\mathcal{P}.
\end{align*}
For instance, in $\overline{\mathbb{H}}^6$, the polytope $\mathcal{P}$ and its dissected simplices are sketched in Fig.~\ref{HoroPolihedron3}.left and \ref{Dissect_Polihedron3} respectively.
\begin{figure}[h!]
\centering 
\includegraphics[width=11cm]{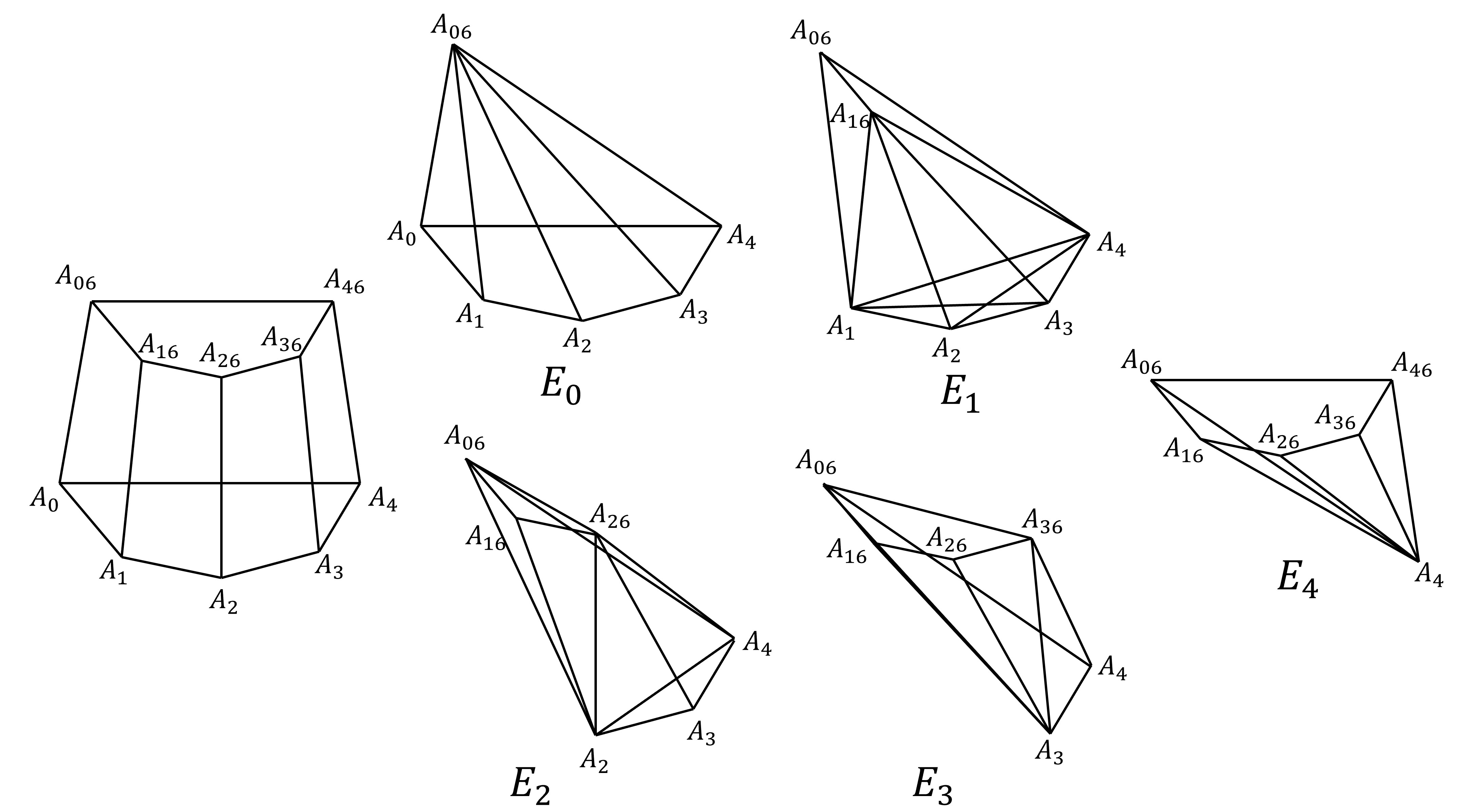}  
\caption{Sketch of the dissection of $5$-dimensional hyperbolical polytope $\mathcal{P}$ into $5$-dimensional simplices $E_k$, $k=0,\cdots 4$.}
\label{Dissect_Polihedron3}
\end{figure}
Furthermore, we take the point projection of each simplex $E_k$ onto the horosphere $\partial\mathcal{B}_{n-1}$ by taking the intersection of $\partial\mathcal{B}_{n-1}$ and the geodesic emanating from the ideal vertex $A_{n-1}$ to all vertices of $\mathcal{P}$. We denote the projection of $E_k$ by $D_k$, respectively. 
To calculate the volumes of simplexes $D_0 \dots D_{n-2}$, 
we use the Cayley-Menger determinant that gives the volume 
of a horospherical $(n-1)$-simplex $\mathrm{vol}(\mathcal{A})$ 
with edges length $L_{i,j}$ ($i=1\dots(n-1), j=2\dots n, i<j$).
\begin{equation}
\mathrm{vol}(\mathcal{A})^2 = \frac{1}{(n!)^22^n}
\begin{vmatrix}\label{Heron}
 0 & 1 & 1 & 1 & \dots &1 \\
 1 & 0 & L_{1,2}^2 & L_{1,3}^2 & \dots & L_{1,n}^2 \\
 1 & L_{1,2}^2 & 0 & L_{2,3}^2 & \dots & L_{2,n}^2 \\
 \vdots & \vdots & \vdots & \ddots & \dots & \vdots \\
 1 & L_{1,n}^2 & L_{2,n}^2 &  \dots & L_{n-1,n}^2 & 0
 \end{vmatrix}.
 \end{equation}

The volume of the horoball piece contained in the fundamental simplex is

\begin{equation}\label{Vol_Hor}
\mathrm{vol}(\mathcal{B}_{n-1} \cap \cF_{\Gamma}) = \frac{1}{n-1}\mathrm{vol}(\mathcal{Q}),~~\text{where}~ \mathcal{Q}:= \bigcup_{k=0\cdots n-2}D_k
\end{equation}
The locally optimal horoball packing density of 
Coxeter Simplex $\cF_{\Gamma}$ can be determined by formula 
$\delta_{opt}(\cF_{\Gamma})  = 
\frac{\mathrm{vol}(\mathcal{B}_{n-1} \cap \cF_{\Gamma})}{\mathrm{vol}(\cF_{\Gamma})}$.

The above local construction is preserved by the isometric actions of $g \in \Gamma$. 
The Coxeter group $\Gamma$ extends the optimal local horoball packing 
density from the fundamental domain $\cF_{\Gamma}$ to the entire 
tiling $\cT_{\Gamma}$ of $\overline{\mathbb{H}}^n$ therefore we obtain the 
horoball packing density of the considered tiling by formulas (4.1) and (4.2)).

\begin{rem}
In some cases, there may be an additional ideal vertex $A_0$, allowing us to consider a horosphere centered at this vertex. The process of determining the density of the optimal horoball packing in these cases is similar to the previously described case. Furthermore, the $(n-1)$-dimensional horospherical polytope $\mathcal{Q}$ is a simplex in these cases, and its volume $\mathrm{vol}(\mathcal{Q})$ can be calculated directly using the Cayley-Menger determinant. 
\end{rem}

\item {\bf Two horoball types}
The computational method for determining the optimal horoball sectors and their volumes is similar to the "one horoball type" method described above. However, in the investigated cases, the two horoballs cannot have a common interior point. At most, they can be tangent to each other. The volumes of the two tangent horoball pieces, centered at two distinct ideal vertices of the fundamental domain, as the horoball type is continuously varied, are related in Lemma \ref{lemma:szirmai}.

In $\overline{\mathbb{H}}^n$ with $n \geq 2$ let $\tau_1$ and $\tau_2$ be two congruent
$n$-dimensional convex cones with vertices at $C_1, C_2 \in \partial \overline{\mathbb{H}}^n$ that share a common geodesic edge $\overline{C_1C_2}$.
Let $B_1(x)$ and $B_2(x)$ denote two horoballs centered at $C_1$ and
$C_2$ respectively, mutually tangent at
$I(x)\in {\overline{C_1C_2}}$. Define $I(0)$ as the point with $V(0) = 2 \mathrm{vol}(B_1(0) \cap \tau_1) = 2 \mathrm{vol}(B_2(0) \cap \tau_2)$ for the volumes of the horoball sectors.

\begin{lem}[\cite{Sz12}]
\label{lemma:szirmai}
Let $x$ be the hyperbolic distance between $I(0)$ and $I(x)$,
then
\begin{equation}
\begin{split}
V(x) =& \mathrm{vol}(B_1(x) \cap \tau_1) + \mathrm{vol}(B_2(x) \cap \tau_2) \\ =& V(0)\frac{e^{(n-1)x}+e^{-(n-1)x}}{2}
= V(0)\cosh\left((n-1)x\right) 
\end{split}
\end{equation}
is strictly convex and strictly increasing as $x\rightarrow\pm\infty$.
\end{lem}
See our paper \cite{Sz12} for a proof. 

The methods described above can be used to determine the optimal horoball packings and their densities.
\end{enumerate}
\subsection{Results related to horoball packings}
\subsubsection{Horoball packings in $\overline{\mathbb{H}}^4$}
\begin{enumerate}
\item {\bf Packings with one horoball}

In the given scenario, we consider a horosphere centered at $A_3$. As mentioned, the horosphere packing reaches its optimal density when the horosphere is tangent to the opposite faces, which are determined by the vertices $A_0, A_1, A_2, A_{04}, A_{14}, A_{24}$. This optimal packing configuration is realized in all possible tilings, as indicated by the Schl\"{a}fli symbols in Fig. 1.

In the case of the Schl\"{a}fli symbol ${4,4,3,4,\infty}$, there are two ideal vertices, namely $A_0$ and $A_3$. Therefore, we can construct another horosphere centered at vertex $A_0$. This horosphere packing will be optimal when it is tangent to the opposite face of $A_0$, which is determined by the vertices $A_1, A_2, A_3, A_4$.

By applying the previously described method, as well as formulas (4.3), (4.4), and (4.5), and utilizing the volume values described in Table 2, we have summarized the results in Table 6.
\begin{table}[h!]
\centering
\begin{footnotesize}
\begin{tabular}{||c c c c||}
		\hline\xrowht[()]{10pt}
		Schl\"afli symbol & $\mathrm{vol}(\mathcal{B}_i \cap \widehat{\mathcal{S}})$ & $\mathrm{vol}(\widehat{\mathcal{S}})$ &  $\delta_{opt}(\mathcal{B}_\Gamma)$.  \\ 
		\hline\hline\xrowht[()]{10pt}
		$ \{4,4,3,3,\infty\},~i=3$ & $ \frac{1}{48}$
		& $ \frac{\pi^2}{288}$  & $\frac{6}{\pi^2}\approx0.60793$
		\\ 
		\hline\xrowht[()]{10pt}
		$\{4,4,3,4,\infty),~i=3$ & $\frac{\sqrt{2}}{48}$
		& $\frac{\pi^2}{144}$ & $\frac{3}{8\pi^2}\approx 0.42987$
		\\
		\hline\xrowht[()]{10pt}
		$\{4,4,3,4,\infty\},~i=0$ & $\frac{\sqrt{2}}{72}$ & $\frac{\pi^2}{144}$ & $\frac{2\sqrt{2}}{\pi^2}\approx {0.28658}$ \\
		\hline\xrowht[()]{10pt}
		$(6,3,3,3,\infty\},~i=3$ & $\frac{\sqrt{3}}{144}$
		& $\frac{\pi^2}{540}$ & $\frac{15\sqrt{3}}{4\pi^2}\approx \textbf{0.65810}$
		\\
		\hline\xrowht[()]{10pt}
		$\{6,3,3,4,\infty\},~i=3$ & $\frac{\sqrt{6}}{144}$
		& $\frac{\pi^2}{288}$ & $\frac{2\sqrt{6}}{\pi^2}\approx 0.49637$
		\\
		\hline\xrowht[()]{10pt}
		$\{6,3,3,5,\infty\},~i=3$ & $\frac{1}{48}\sqrt{\frac{3+\sqrt{5}}{6}}$
		& $\frac{61\pi^2}{900}$ & $\frac{25(\sqrt{3}+\sqrt{15})}{244\pi^2}\approx 0.02909$
		\\
		\hline\xrowht[()]{10pt}
		$\{6,3,4,3,\infty\},~i=3$ & $\frac{\sqrt{3}}{72}$
		& $\frac{5\pi^2}{864}$ & $\frac{12\sqrt{3}}{5 \pi^2}\approx 0.42119$
		\\
		\hline
	\end{tabular}
\end{footnotesize}
 \caption{Horoball Packing densities in $\overline{\mathbb{H}}^4$, with one horosphere}
 \label{densH4_1}
\end{table}
\newpage
\item {\bf Packings with two horoballs}
\begin{thm}
The optimal horoball packing density of Coxeter simplex tilings
$\cT_\Gamma$ where $\Gamma$ is determined by  
Schl\"{a}fli symbol $\{4,4,3,4,\infty\}$ is 
\begin{equation}
\delta_{opt}(\cT_\Gamma) = 
\frac{\mathrm{vol}(\mathcal{B}_0^{opt} \cap \widehat{\mathcal{S}})+
\mathrm{vol}(\mathcal{B}_3^{opt} \cap \widehat{\mathcal{S}}))}{\frac{\pi^2}{144}} 
\approx 0.71645.
\end{equation}
\end{thm}
\textit{Proof} \\
We consider the case of Schl\"{a}fli symbol $\{4,4,3,4,\infty\}$ whose the Coxeter-Schl\"{a}fli matrix is the following Gram matrix
\begin{equation}\label{Gram4}
    \displaystyle (b^{ij})=(\langle \boldsymbol{b}^i, \boldsymbol{b}^j \rangle )= \left[\begin{matrix}1 & - \frac{\sqrt{2}}{2} & 0 & 0 & 0\\- \frac{\sqrt{2}}{2} & 1 & - \frac{\sqrt{2}}{2} & 0 & 0\\0 & - \frac{\sqrt{2}}{2} & 1 & - \frac{1}{2} & 0\\0 & 0 & - \frac{1}{2} & 1 & - \frac{\sqrt{2}}{2}\\0 & 0 & 0 & - \frac{\sqrt{2}}{2} & 1\end{matrix}\right],
\end{equation}
where $0 \leq i,j \leq 4$.
Moreover, the corresponding inverse matrix $(a_{ij})$ 
\begin{equation}
    \displaystyle (a_{ij})= (\langle \boldsymbol{b}^i, \boldsymbol{b}^j \rangle )^{-1}= \left[\begin{matrix}0 & - \sqrt{2} & -2 & -2 & - \sqrt{2}\\- \sqrt{2} & -2 
    & - 2 \sqrt{2} & - 2 \sqrt{2} & -2\\-2 & - 2 \sqrt{2} & -2 & -2 & - \sqrt{2}\\-2 & - 2 \sqrt{2} & -2 & 0 & 0\\- \sqrt{2} & -2 & - \sqrt{2} & 0 & 1\end{matrix}\right]
\end{equation}
shows that the truncated orthoscheme has two ideal vertices $A_0[\bold{a}_0]$, and $A_3[\bold{a}_3]$  since $a_{00}=0$, and $a_{33}=0$. We have an ultra ideal vertices  $A_4[\bold{a}_4]$, since $a_{44}= 1 > 0$.\\
First, we choose the set unit normal forms of the orthoscheme hyperplane faces 
$\boldsymbol{b}^0=(0,0,1,0,0)^T$, $\boldsymbol{b}^1=(-\frac{\sqrt{2}}{2},,-\frac{\sqrt{2}}{2},-\frac{\sqrt{2}}{2},-\frac{\sqrt{2}}{2},0)^T$, $\boldsymbol{b}^2=(0,0,0,1,0)^T$, 
$\boldsymbol{b}^3=(0,\frac{1}{2},0,-\frac{1}{2},-\frac{\sqrt{2}}{2})^T$, $\boldsymbol{b}^4=(0,0,0,0,1)^T$ admitted the Gram matrix (Coxeter Schl\"{a}fli Matrix)  (\ref{Gram4}).\\
Those unit normal vectors $\{\boldsymbol{b}^0,\boldsymbol{b}^1,\boldsymbol{b}^2,\boldsymbol{b}^3,\boldsymbol{b}^4\}$ and $(a_{ij})$ determine the orthoscheme vertices whose coordinates: $\bold{a}_0=(1,0,1,0,0)$, $\bold{a}_1=(1,0,0,0,0)$, $\bold{a}_2=(1,\frac{1}{2},0,\frac{1}{2},0)$, $\bold{a}_3=(1,1,0,0,0)$, $\bold{a}_4=(1,1,0,0,\frac{\sqrt{2}}{2})$.
We consider two horoballs $\mathcal{B}_0$ and $\mathcal{B}_3$ centred at $A_0$ and $A_3$ which hold the packing conditions, so they may touch each other 
at a point on edge $A_0 A_3$. 
To investigate this possibility, we consider a linear parametrization of moving point $P[\mathbf{p}]$ points on $A_0 A_3$
\begin{equation}\label{eq:P_30}
    \bold{p}(t)=(1-t)\bold{a}_3 + t \bold{a}_0, \qquad t \in [0,1],
\end{equation}
where clearly $\bold{p}(0)=\bold{a}_3$ and $\bold{p}(1)=\bold{a}_0$.  
Based on the method described in subsection 4.1 we can determine possible maximal horospheres centred at $A_0$ and 
$A_3$ and its parameters $s_i, i\in\{0,3\}$ using the projection formula (4.3) and the equation of horospheres given in (2.4). 
 \begin{equation*}
     \mathbf{f}_3 =\ba_{3} - \frac{\langle \ba_{3}, \bb^3 \rangle}{\langle \Bb^3,\Bb^3 \rangle} 
\bb^3=\left(1,\frac{3}{4},0,\frac{1}{4}, \frac{\sqrt{2}}{4}\right)
 \end{equation*}
 An isometry could take $A_3$ to the point $(1,0,0,0,1)$, and $\bold{f}_3$ to $\hat{\bold{f}}_3=(1,\frac{\sqrt{2}}{4},\frac{1}{4},0,\frac{3}{4})$
 by substituting $\hat{\bold{f}}_3$ to the horosphere equation given in (\ref{eq:horosphere}), we have
 \begin{equation*}
     \frac{\left(\frac{1}{2}-s\right)^2}{\left(1-s \right)^2}+\frac{\frac{3}{8}}{1-s}=1
 \end{equation*}
 by solving the last equation for $s$, then it is concluded that $s_3=\frac{3}{5}$, furthermore the corresponding $t$ in (\ref{eq:P_30}) is $t_3=\frac{1}{3}$, i.e., $\bold{p}(\frac{1}{3})$ is the intersection of $\mathcal{B}_3$ and the edge $A_3A_0$.\\
On the other hand, We consider the horoball $\mathcal{B}_0$ centred at $A_0$  such that it tangent either the opposite face, $H_0$ determined by vertices $A_1,A_2,A_3,A_4$ or the truncation face $H_5$ given by vertices $A_{04}, A_{14}, A_{24}, A_3$. We consider the following projections:
\begin{equation*}
    \mathbf{f}_0 =\ba_{0} - \frac{\langle \ba_{0}, \bb^0 \rangle}{\langle \Bb^0,\Bb^0 \rangle} \bb^0=(1,0,0,0,0)=\ba_1
\end{equation*}
There is an isometry that takes $\ba_0$ to $(1,0,0,0,1)$, where $\bold{f}_0=\ba_1$ is the fix point. By substituting $\bold{f}_0$ into the equation (\ref{eq:horosphere}), and solving the equation for $s$ we obtain $s_0=0$. This parameter $s_0=0$ determine the horosphere $\mathcal{B}_0^{opt}$.\\
Consider 
\begin{equation*}
    \bold{q}(t)=(1-t)\bold{a}_0+t\bold{a}_3,~~t \in [0,1]
\end{equation*}
this is obviously the opposite of $\bold{p}$, i.e., $\bold{q}(t)=\bold{p}(1-t)$. By careful analysis, the intersection of $\mathcal{B}_0(s_0)$ and the edge $A_0A_3$ is the point $\bold{q}(\frac{2}{3})$. As a result,
\begin{equation*}
    \bold{q}\left(\frac{2}{3}\right)=\bold{p}\left(1-\frac{2}{3}\right)=\bold{p}\left(\frac{1}{3}\right)
\end{equation*}
that means the horoballs $\mathcal{B}_0^{opt}$ and $\mathcal{B}_3^{opt}$ 
touch each other at a point on edge $A_0A_3$.

The optimum density of horoball packing with two horospheres can be directly computed by (\ref{Heron}) and (\ref{Vol_Hor})
\begin{equation}
\begin{gathered}
\delta_{opt}(\cT_\Gamma) = \frac{\mathrm{vol}(\mathcal{B}_0^{opt} \cap \widehat{\mathcal{S}})+\mathrm{vol}(\mathcal{B}_3^{opt} 
\cap \widehat{\mathcal{S}}))}{\frac{\pi^2}{144}} \\
= \frac{\mathrm{vol}(\mathcal{B}_0(0) \cap \widehat{\mathcal{S}})+\mathrm{vol}(\mathcal{B}_3(\frac{3}{5}) 
\cap \widehat{\mathcal{S}}))}{\frac{\pi^2}{144}}
=\frac{\frac{\sqrt{2}}{48}+\frac{\sqrt{2}}{72}}{\frac{\pi^2}{144}}=\frac{5\sqrt{2}}{\pi^2}\approx 0.71645~~~\square
\end{gathered} \notag
\end{equation}
Finally, we summarize our results in the following theorem
\begin{thm}
In hyperbolic space $\overline{\mathbb{H}}^4$, between the congruent ball and ho\-ro\-ball packings of at most two ho\-ro\-ball types, 
generated by simply truncated Coxeter orthoschemes with parallel faces, the above determined $\mathcal{B}_{\{4,4,3,4,\infty\}}$ horoball configuration 
provides the densest packing with density $\frac{5\sqrt{2}}{\pi^2}\approx 0.71645$. 
\end{thm}
\begin{rem}
This density $\delta_{opt}(\cT_{\Gamma}) \approx 0.71645$ is the known densest 
horoball packing density in $\overline{\mathbb{H}}^4$.  
We note here, that this density is realized in seven asymptotic Coxeter 
simplex tilings $\Gamma \in \{ \overline{S}_4, \overline{P}_4, \widehat{FR}_4, 
\overline{N}_4, \overline{M}_4, \overline{BP}_4,$ $\overline{DP}_4 \}$, 
if horoballs of different types are allowed at each asymptotic vertex 
of the tiling (see \cite{KSz14}).
\end{rem}
\end{enumerate}
\subsubsection{Horoball packings in $\overline{\mathbb{H}}^5$}
The Coxeter Schl\"{a}fli matrix of orthoscheme whose Schl\"{a}fli symbol 
$\{4,3,4,3,3,\infty\}$ is $(b^{ij})$. 
Choosing the unit normal forms of hyperplane faces of the orthoscheme admitted 
the matrix $(b^{ij})$, i.e.,
\begin{equation}
\begin{gathered}
\Bb^0=(0,0,1,0,0,0)^T, \Bb^1=\left(-\frac{\sqrt{15}}{3}-\sqrt{3}, -\frac{3\sqrt{2}}{2}-\frac{\sqrt{10}}{2},-\frac{\sqrt{2}}{2},-\frac{1}{2},0,0\right)^T,\\
\Bb^2=(0,0,0,1,0,0)^T, \Bb^3=\left(\frac{\sqrt{30}}{6},1,0,-\frac{\sqrt{2}}{2},-\frac{\sqrt{3}}{3},0\right)^T, \\
\Bb^4=\left( 0,0,0,0,\frac{\sqrt{3}}{2},-\frac{1}{2}\right)^T, \Bb^5=(0,0,0,0,0,1)^T.
\end{gathered}
\end{equation}
The corresponding vertices are given by $A_0[\Ba_0]$, $A_1[\Ba_1],$
$A_2[\Ba_2]$, $A_3[\Ba_3], A_4[\Ba_4],$ $A_5[\Ba_5]$, where 
\begin{equation}
\begin{gathered}
\ba_0=\left(1, \frac{\sqrt{30}}{6}, \frac{\sqrt{6}}{6},0,0,0\right), \ba_1=\left( 1, \frac{\sqrt{30}}{6},0,0,0,0 \right),\\
\ba_2=\left(1, \frac{2\sqrt{6}+5\sqrt{30}}{33},0, \frac{4\sqrt{3}-\sqrt{15}}{33},0,0 \right), \ba_3=\left(1, \frac{\sqrt{6}+\sqrt{30}}{8},0,0,0,0 \right)\\
\ba_4=\left(1,\frac{\sqrt{6}+\sqrt{30}}{8},0,0,\frac{3\sqrt{2}-\sqrt{10}}{8},0\right), \\
\ba_5=\left(1,\frac{\sqrt{6}+\sqrt{30}}{8},0,0,\frac{3\sqrt{2}-\sqrt{10}}{8},\frac{3\sqrt{6}-\sqrt{30}}{8}  \right).
\end{gathered}
\end{equation}
\begin{enumerate}
\item {\bf Packings with one horoball}

Figure 2 illustrates that in the $5$-dimensional hyperbolic space $\overline{\mathbb{H}}^5$, there is only one Coxeter tiling $\mathcal{T}_\Gamma$ generated by a simply truncated orthoscheme with parallel faces of Schl\"afli symbol ${4,3,4,3,3,\infty}$.
By observing the corresponding $(a{ij})=(b^{ij})^{-1}$ matrix, we note that its fundamental domain $\mathcal{F}_\Gamma$ has two ideal vertices, $A_0$ and $A_4$. Hence, we can construct horoballs centered at either of these ideal vertices. In the case of the horosphere centered at $A_0$, the opposite face is $A_1 A_2 A_3 A_4 A_5$, and the truncation face is $A_4 A{05} A_{15} A_{25} A_{35}$. However, the optimal packing configuration can only be realized when the horosphere touches the face $A_1 A_2 A_3 A_4 A_5$.
Similar to the $4$-dimensional case described above, by applying formulas (4.3), (4.4), (4.5), and the volume value determined in subsection 2.3.2, we can determine the densest horoball packings with one horoball type. The results are summarized in Table 7.
\begin{table}[h!]
	\centering
	\begin{footnotesize}
		\begin{tabular}{||c c c c||}
			\hline\xrowht[()]{10pt}
			Schl\"afli symbol & $\mathrm{vol}(\mathcal{B}_i \cap \widehat{\mathcal{S}})$ & $\mathrm{vol}(\widehat{\mathcal{S}})$ &  $\delta_{opt}(\mathcal{B}_\Gamma)$.  \\ 
			\hline\hline\xrowht[()]{10pt}
			$\{4,3,4,3,3,\infty\},~i=4$ & $\frac{1}{384}$
			
			& {$\frac{7\zeta{(3)}}{1152}$} & $ \frac{3}{7\zeta{(3)}}\approx \bold{0.35653}$
			\\ 
			\hline\xrowht[()]{10pt}
			$\{4,3,4,3,3,\infty\},~i=0$ & $\frac{1}{576}$
			&{$\frac{7\zeta{(3)}}{1152}$} & $ \frac{2}{7\zeta{(3)}}\approx 0.23769$
			\\
			\hline
		\end{tabular}
	\end{footnotesize}
	 \caption{Packing densities in $\overline{\mathbb{H}}^5$, with one horoball}
	\label{densH5_1}
\end{table}
\newpage
\item {\bf Packings with two horoballs}

Here we consider two horospheres centered at $A_0$ and $A_4$ together and compute their optimal arrangements and density. It is analogous to the previous construction in 
$4$-dimensional hyperbolic space (see 4.2.1 subsection). We define the following linear parametrization
\begin{equation}
	\bold{p}(t)=(1-t)\bold{a}_4 + t \bold{a}_0, \ \ t \in [0,1],
\end{equation}
that describes points along edge $A_0 A_4$. By careful investigation, we find that if $t= \frac{27-6\sqrt{5}}{61}$, then the horosphere centred at $A_4$ osculates the opposite face.
Furthermore, this exact value of $t$ provides the conditions of the other horosphere centred at $A_0$ osculate the opposite face $A_1 A_2 A_3 A_4 A_5$. 
Therefore, these two horospheres tangent each other if and only if they osculate their centres opposite faces. We summarize the computation result in Table 8.
\begin{table}
	\centering
	\begin{footnotesize}
		\begin{tabular}{||c c c c||}
			\hline\xrowht[()]{10pt}
			Schl\"afli symbol & $\mathrm{vol}((\mathcal{B}_0 \cup \mathcal{B}_4) \cap \widehat{\mathcal{S}})$ & $\mathrm{vol}(\widehat{\mathcal{S}})$ &  
			$\delta_{opt}(\mathcal{B}_\Gamma))$.  \\ 
			\hline\hline\xrowht[()]{10pt}
			$(4,3,4,3,3,\infty),~i=0,4$ & $\frac{5}{1152}$
			
			& {$\frac{7\zeta{(3)}}{1152}$} & $\frac{5}{7\zeta{(3)}} \approx \bold{0.59421}$
			\\ 
			
			\hline
		\end{tabular}
	\end{footnotesize}
	 \caption{Packing densities in $\overline{\mathbb{H}}^5$, with two horoballs}
	\label{densH5_2}
\end{table}
We summarize our investigation in 5-dimensional case with the following theorem.
\begin{thm}
In hyperbolic space $\overline{\mathbb{H}}^5$, between the congruent ball and ho\-ro\-ball packings of at most two ho\-ro\-ball types, generated by simply 
truncated Coxeter orthoschemes with parallel faces, 
the above determined $\mathcal{B}_{\{4,3,4,3,3,\infty\}}$ horoball configuration 
provides the densest packing with density $\frac{5}{7\zeta{(3)}}\approx 0.59421$.
\end{thm}
\end{enumerate}
\begin{rem}
In \cite{KSz18} we proved that the above optimal density is realized in 
ten commensurable asymptotic Coxeter simplex tilings 
$$\Gamma \in \left\{\overline{U}_5,  \overline{S}_5, \overline{Q}_5, \overline{X}_5, \overline{R}_5, \overline{N}_5, \overline{O}_5, \overline{M}_5, \overline{L}_5, \widehat{UR}_5  \right\},$$ 
when horoballs of different types are allowed at each asymptotic vertex of the tiling. 
\end{rem}
\newpage
\subsubsection{Horoball packings in $\overline{\mathbb{H}}^6$}
\begin{enumerate}
\item {\bf Packings with one horoball}

Tha Fig.~2 shows that in the $6$-dimensional hyperbolic space $\overline{\mathbb{H}}^6$ there are two Coxeter tiling $\mathcal{T}_{\Gamma_i}, ~ i\in\{1,2\}$ generated by 
simply truncated orthosche\-mes with parallel faces of Sch\"afli symbols 
$\{3,4,3,3,3,3,\infty\}$ and $\{3,4,3,3,3,4,$ $\infty\}$.
Using the corresponding $(a_{ij})=(b^{ij})^{-1}$ matrices it is clear, that in the first case the fundamental domain $\mathcal{F}_{\Gamma_1}$ has only one ideal vertex $A_5$ 
and in the second case the fundamental domain $\mathcal{F}_{\Gamma_1}$ has has two ideal vertices, 
$A_0$ and $A_5$. Hence, we can construct horoballs centred at either 
ideal vertices.  
The horosphere should be created, centred at ideal vertex $A_5$, as mentioned in the previous subsections. 
Furthermore, in configuration of Schl\"{a}fli symbol $\{3,4,3,3,3,4,\infty\}$, we have another horosphere centred at $A_0$.
The horoball sector centred at $A_0$ will reach maximum volume, if it osculates the opposite face of $A_0$, i.e., the face $A_1 A_2 A_3 A_4 A_5 A_6$. 
On the other hand, the hull of vertices $A_1 A_2 A_3 A_4 A_5 A_6$ is an 5-simplex in $\overline{\mathbb{H}}^6$. We can take the central projection of this simplex on the horosphere, 
and compute the volume of its horoball sector. Similarly to the above $4$- and $5$-dimensional case applying the formulas 
(4.3), (4.4), (4.5), and the volume values in Table 2. We can determine the densest horoball packings with one horoball type. We summarized the results in Table 9.                         
\begin{table}
\centering
\begin{footnotesize}
	\begin{tabular}{||c c c c||}
		\hline\xrowht[()]{10pt}
		Schl\"afli symbol & $\mathrm{vol}(\mathcal{B}_i \cap \hat{\mathcal{S}})$ & $\mathrm{vol}(\hat{\mathcal{S}})$ &  $\delta_{opt}(\mathcal{B}_\Gamma))$.  \\ 
		\hline\hline\xrowht[()]{10pt}
		$(3,4,3,3,3,3,\infty),~i=5$ & $\frac{1}{23040}$ 
		
		& $\frac{11\pi^3}{86400}$ & $\approx 0.01010$
		\\ 
		\hline\xrowht[()]{10pt}
		$(3,4,3,3,3,4,\infty),~i=5$ & $\frac{\sqrt{2}}{23040}$
		& $\frac{\pi^3}{86400}$ & $\approx \bold{0.17114}$
		\\
		\hline\xrowht[()]{10pt}
		$(3,4,3,3,3,4,\infty),~i=0$ & $\frac{1}{19200}$
		& $\frac{\pi^3}{86400}$ & $\approx 0.14513$
		\\
		\hline
	\end{tabular}
\end{footnotesize}
	 \caption{Packing densities in $\overline{\mathbb{H}}^6$, with one horoballs}
	\label{densH6_1}
\end{table}
\item {\bf Packings with two horoballs}

In case of the tiling $\mathcal{T}_{\Gamma_2}$ given by Schl\"{a}fli symbol $\{3,4,3,3,3,4,\infty\}$, we have two ideal vertices $A_0$ and $A_5$. 
These two horospheres might tangent each other. The point of tangency lies on the edge $A_0 A_5$ if it exists.
We can determine the maximum possible horoballs in the fundamental domain (truncated orthoscheme) $\mathcal{F}_{\Gamma_2}$ belonging to the ideal vertices, 
we find that they have no common point, they are disjoint. Therefore the optimal horoball arrangement is derived from this horoball configuration.
The optimum density of horoball packings with two horoballs can be directly computed by
	\begin{equation}
	\begin{gathered}
	\delta_{opt}(\cT_{\Gamma_2}) = \frac{\mathrm{vol}(\mathcal{B}_0^{opt} \cap \widehat{\mathcal{S}})+\mathrm{vol}(\mathcal{B}_5^{opt} 
	\cap \widehat{\mathcal{S}}))}{\frac{ \pi^3}{86,400}}=\\
    =\frac{\mathrm{vol}(\mathcal{B}_0(0) \cap \widehat{\mathcal{S}})+\mathrm{vol}(\mathcal{B}_5(\frac{7}{9}) 
	\cap \widehat{\mathcal{S}}))}{\frac{ \pi^3}{86,400}}=\\
    =\frac{\frac{\sqrt{2}}{23040}+\frac{1}{19200}}{\frac{ \pi^3}{86,400}}=\frac{15\sqrt{2}}{4\pi^3}+\frac{9}{2\pi^3} 
	\approx 0.31617 ~~~\square
	\end{gathered} \notag
\end{equation}
We summarize our results in the following table and theorem.
\begin{table}
	\centering
\begin{footnotesize}
	\begin{tabular}{||c c c c||}
			\hline\xrowht[()]{10pt}
			Schl\"afli symbol & $\mathrm{vol}((\mathcal{B}_0 \cup \mathcal{B}_5) \cap \hat{\mathcal{S}})$ & 
			$\mathrm{vol}(\hat{\mathcal{S}})$ &  $\delta_{opt}(\mathcal{B}_\Gamma))$.  \\ 
			\hline \hline\xrowht[()]{10pt}
			$\{3,4,3,3,3,4,\infty\}$ & $\frac{5+3\sqrt{2}}{57600\sqrt{2}}$ & $\frac{\pi^3}{86400}$ & $\approx \bold{0.31617}$ \\ 
		\hline
		\end{tabular}
	\end{footnotesize}
	 \caption{Packing densities in $\overline{\mathbb{H}}^6$, with two horospheres}
	\label{densH6_2}
\end{table}
\begin{thm}
	In hyperbolic space $\overline{\mathbb{H}}^6$, between the congruent ball and ho\-ro\-ball packings of at most two ho\-ro\-ball types, generated by simply 
	truncated Coxeter orthoschemes with parallel faces, 
	the above determined $\mathcal{B}_{\{3,4,3,3,3,4, \infty \}}$ horoball configuration 
provides the densest packing with density $\frac{15\sqrt{2}+18}{4\pi^3}\approx 0.31617$.
\end{thm}
\end{enumerate}
The ball, horoball, hyperball packings, and coverings in Thurston geometries as well as in higher dimensional hyperbolic spaces still contain many open questions, which we plan to investigate further, which may lead to further interesting results.
%

\end{document}